\documentclass{amsart}

\let\hat=\widehat

\newcommand{\C}{\mathbb{C}}
\newcommand{\F}{\mathbb{F}}
\newcommand{\Q}{\mathbb{Q}}
\newcommand{\Z}{\mathbb{Z}}

\def\calO{\mathcal{O}}
\def\calL{\mathcal{L}}

\def\frakp{\mathfrak{p}}

\newcommand{\alphabar}{\overline{\alpha}}
\newcommand{\betabar}{\overline{\beta}}
\newcommand{\pibar}{\overline{\pi}}
\newcommand{\rhobar}{\overline{\rho}}
\newcommand{\qbar}{\overline{q}} 
\newcommand{\rbar}{\overline{r}} 
\newcommand{\Fpbar}{\overline{\F}_p}
\newcommand{\Qb}{\overline{\Q}}
\newcommand{\kbar}{\overline{k}}

\newcommand{\psihat}{\hat{\psi}}
\newcommand{\Kp}{K^+} 
\newcommand{\Grr}{G_{\rm{red,red}}}
\newcommand{\Grl}{G_{\rm{red,loc}}}
\newcommand{\Glr}{G_{\rm{loc,red}}}
\newcommand{\Gll}{G_{\rm{loc,loc}}}

\DeclareMathOperator{\divisor}{div}
\DeclareMathOperator{\End}{End}
\DeclareMathOperator{\Hom}{Hom}
\DeclareMathOperator{\Jac}{Jac} 
\DeclareMathOperator{\Norm}{Norm} 
\DeclareMathOperator{\ord}{ord} 
\DeclareMathOperator{\trace}{trace} 

\newcommand\Bdownarrow{\vbox{\vbox to 4pt{}\vbox{\Big\downarrow\vfill}}}
\newcommand\Buparrow  {\vbox{\vbox to 4pt{}\vbox{\Big\uparrow  \vfill}}} 
\newcommand\mapright[1]{\mathop{\longrightarrow}\limits^{#1}}

\newtheorem{thm} {Theorem}
\newtheorem{lem} [thm]{Lemma} 
\newtheorem{cor} [thm]{Corollary}
 \newtheorem{prop}[thm]{Proposition}

\newtheorem*{claim}{Claim}
\newtheorem*{assumption}{Assumption}

\begin{document}

\title[Improved upper bounds for curves]
      {Improved upper bounds for the number\\ 
       of points on curves over finite fields}

\author{Everett W.~Howe}
\address{Center for Communications Research, 
         4320 Westerra Court, 
         San Diego, CA 92121-1967, USA.} 
\email{however@alumni.caltech.edu}
\urladdr{http://alumni.caltech.edu/\~{}however/}

\author{Kristin E. Lauter}
\address{Microsoft Research,
         One Microsoft Way,
         Redmond, WA 98052, USA.}
\email{klauter@microsoft.com}

\date{17 December 2006}



\begin{abstract}
We give new arguments that improve the known upper bounds on 
the maximal number $N_q(g)$ of rational points of a curve of 
genus $g$ over a finite field~$\F_q$, for a number of pairs $(q,g)$. 
Given a pair $(q,g)$ and an integer~$N$, we determine the possible 
zeta functions of genus-$g$ curves over $\F_q$ with $N$ points,
and then deduce properties of the curves from their zeta functions.
In many cases we can show that a genus-$g$ curve over $\F_q$ with 
$N$ points must have a low-degree map to another curve over~$\F_q$,
and often this is enough to give us a contradiction.
In particular, we are able to provide eight previously unknown
values of $N_q(g)$, namely:  
$N_4(5) = 17$, $N_4(10) = 27$, $N_8(9) = 45$, $N_{16}(4) = 45$, 
$N_{128}(4) = 215$, $N_3(6) = 14$, $N_9(10) = 54$, and $N_{27}(4) = 64$.
Our arguments also allow us to give a non-computer-intensive proof
of the recent result of Savitt that there are no genus-$4$ curves 
over $\F_8$ having exactly $27$ rational points. 
Furthermore, we show that there is an infinite sequence of 
$q$'s such that for every $g$ with $0<g<\log_2 q$,
the difference between the Weil-Serre
bound on $N_q(g)$ and the actual value of $N_q(g)$ is at least $g/2$.
\end{abstract}

\maketitle

\section{Introduction}
\label{S-intro}
The number $N$ of points on a smooth, absolutely irreducible curve 
of genus $g$ over a finite field $\F_q$ is bounded by 
$$q+1-2g\sqrt{q} \le N \le q+1+2g\sqrt{q},$$ 
an estimate given by Andr\'e Weil in the 1940s.  
In 1983, Serre improved this bound to 
$$q+1-gm \le N \le q+1+gm, \quad \text{where $m=[2\sqrt{q}]$.}$$ 
Serre also introduced the {\em explicit formul{\ae} method}, which uses
numerical conditions on the number of points on a curve over extensions
of the ground field to obtain improved bounds on~$N$, 
at least when $g$ is large compared to $q$ (specifically, 
when $g > (q - \sqrt{q})/2$). 
Oesterl\'e optimized the explicit formul{\ae} method, and 
the resulting bounds on $N$ are the best possible bounds that can 
be obtained formally using only Weil's ``Riemann hypothesis'' for curves
and the fact that for every $d\ge0$ the number of places
of degree $d$ on a curve is non-negative. 
But the method does not take the geometry of the curves into account,
and for this reason it is natural to suspect that the Weil-Serre-Oesterl\'e
bounds may not be optimal.  
Indeed, Serre~\cite{Serre:notes} and others 
\cite{FuhrmannTorres,
KorchmarosTorres, 
Lauter:CRAS,
Lauter:PAMS,
Lauter:conf,
LauterSerre:JAG,
LauterSerre:CM,
Savitt,
Serre:CRAS,
Serre:Bordeaux,
Serre:Resume,
Stark, 
StohrVoloch, 
Zieve}
have shown that in certain cases these bounds are not attained.  
However, in many other cases the bounds provided by the explicit 
formul{\ae} method remain the best known, and significant effort has been made
to determine whether or not they are met --- see the tables of
van der Geer and van der Vlugt~\cite{GeerVlugt}, 
which summarize the work of many authors.  
The purpose of this paper is to provide some new techniques
that show that in many cases the current upper bounds cannot be met.  
In particular, we list in Tables~\ref{Tbl-char2} 
and~\ref{Tbl-char3} the
improvements we obtain to the tables in~\cite{GeerVlugt}. 
In Table~\ref{Tbl-exact} we list the values of $q$ and $g$ for
which our improved upper bound meets the known lower bound.
(The tables in~\cite{GeerVlugt} are updated frequently;
the latest version can be found~at

\smallskip
\noindent
{\tt http://www.science.uva.nl/\~{}geer/}
\smallskip

\noindent
The version of the tables that we will refer to in this paper
is dated 18 January 2002, and is available at the URL mentioned
below in the acknowledgments.)

\begin{table}
\begin{center}
\begin{tabular}{|c|c|c|c|}
\hline
  $q$ & genus            & old upper bound & new upper bound  \\ \hline
  \hline
  $4$ & $5$              & $18$            & $17$             \\ \hline
  $4$ & $10$             & $28$            & $27$             \\ \hline
  $4$ & $11$             & $30$            & $29$             \\ \hline
  \hline
  $8$ & $5$              & $32$            & $30$             \\ \hline
  $8$ & $7$              & $39$            & $38$             \\ \hline
  $8$ & $8$              & $43$            & $42$             \\ \hline
  $8$ & $9$              & $47$            & $45$             \\ \hline
  $8$ & $10$             & $50$            & $49$             \\ \hline
  $8$ & $11$             & $54$            & $53$             \\ \hline
  $8$ & $15$             & $68$            & $67$             \\ \hline
  \hline
 $16$ & $4$              & $46$            & $45$             \\ \hline
 $16$ & $5$              & $54$            & $53$             \\ \hline
 $16$ & $7$              & $70$            & $69$             \\ \hline
 $16$ & $8$              & $76$            & $75$             \\ \hline
 $16$ & $11$             & $92$            & $91$             \\ \hline
 $16$ & $13$             & $103$           & $102$            \\ \hline
 $16$ & $14$             & $108$           & $107$            \\ \hline
  \hline
 $32$ & $4 \le g \le  8$ & $q+1+gm-2$      & $q+1+gm-3$       \\ \hline
 $32$ & $9 \le g \le 15$ & $q+1+gm-2$      & $q+1+gm-4$       \\ \hline
  \hline
 $64$ & $11\le g\le 27$  & $q+1+gm-3$      & $q+1+gm-5$       \\
      & and $g \neq 12$  &                 &                  \\ \hline 
  \hline
$128$ & $4$              & $217$           & $215$            \\ \hline
$128$ & $6$              & $261$           & $258$            \\ \hline
$128$ & $8$              & $305$           & $302$            \\ \hline
$128$ & $9$              & $327$           & $322$            \\ \hline
$128$ & $11$             & $371$           & $366$            \\ \hline
$128$ & $15\le g \le 64$ & $q+1+gm$        & $q+1+gm-4$       \\ 
      & and $g\equiv1\bmod 7$ &            &                  \\ \hline
$128$ & $16\le g \le 65$ & $q+1+gm$        & $q+1+gm-5$       \\ 
      & and $g\equiv2\bmod 7$ &            &                  \\ \hline
$128$ & $10\le g \le 59$ & $q+1+gm$        & $q+1+gm-4$       \\ 
      & and $g\equiv3\bmod 7$ &            &                  \\ \hline
$128$ & $18\le g \le 60$ & $q+1+gm$        & $q+1+gm-6$       \\ 
      & and $g\equiv4\bmod 7$ &            &                  \\ \hline
$128$ & $5 \le g \le 61$ & $q+1+gm$        & $q+1+gm-5$       \\ 
      & and $g\equiv5\bmod 7$ &            &                  \\ \hline
$128$ & $13\le g \le 62$ & $q+1+gm$        & $q+1+gm-7$       \\ 
      & and $g\equiv6\bmod 7$ &            &                  \\ \hline
\end{tabular}
\end{center} 
\vspace{1ex}
\caption{Improved upper bounds on the number of points on curves of 
certain genera over small finite fields~$\F_q$ of characteristic~$2$.
The symbol $m$ is an abbreviation for $[2\sqrt{q}]$.} 
\label{Tbl-char2} 
\end{table}

\begin{table}
\begin{center}
\begin{tabular}{|c|c|c|c|}
\hline
  $q$ & genus            & old upper bound & new upper bound  \\ \hline
  \hline
  $3$ & $6$              & $15$            & $14$             \\ \hline
  \hline
  $9$ & $9$              & $51$            & $50$             \\ \hline
  $9$ & $10$             & $55$            & $54$             \\ \hline
  $9$ & $11$             & $59$            & $58$             \\ \hline
  $9$ & $12$             & $63$            & $62$             \\ \hline
  $9$ & $13$             & $66$            & $65$             \\ \hline
  $9$ & $14$             & $70$            & $69$             \\ \hline
  $9$ & $15$             & $74$            & $73$             \\ \hline
  $9$ & $16$             & $78$            & $77$             \\ \hline
  $9$ & $17$             & $82$            & $81$             \\ \hline
  $9$ & $18$             & $85$            & $84$             \\ \hline
  \hline
 $27$ & $4$              & $66$            & $64$             \\ \hline
 $27$ & $5 \le g \le 8$  & $q+1+gm-2$      & $q+1+gm-3$       \\ \hline
 $27$ & $9 \le g \le 13$ & $q+1+gm-2$      & $q+1+gm-5$       \\ \hline
 $27$ & $14$             & $164$           & $163$            \\ \hline
  \hline
 $81$ & $13 \le g \le 17$& $q+1+gm-2$      & $q+1+gm-4$       \\
      & and $g \neq 16$  &                 &                  \\ \hline
 $81$ & $18\le g \le 35$ & $q+1+gm-2$      & $q+1+gm-5$       \\ \hline
\end{tabular}
\end{center} 
\vspace{1ex}
\caption{Improved upper bounds on the number of points on curves of 
certain genera over small finite fields~$\F_q$ of characteristic~$3$.
The symbol $m$ is an abbreviation for $[2\sqrt{q}]$.} 
\label{Tbl-char3} 
\end{table}

\begin{table}
\begin{center}
\begin{tabular}{|c|c|c|}
\hline
  $q$ & genus $g$ & $N_q(g)$ \\ \hline
  \hline
  $4$ & $ 5$      & $ 17$    \\ \hline
  $4$ & $10$      & $ 27$    \\ \hline
  $8$ & $ 9$      & $ 45$    \\ \hline
 $16$ & $ 4$      & $ 45$    \\ \hline
$128$ & $ 4$      & $215$    \\ \hline
  $3$ & $ 6$      & $ 14$    \\ \hline
  $9$ & $10$      & $ 54$    \\ \hline
 $27$ & $ 4$      & $ 64$    \\ \hline
\end{tabular}
\end{center} 
\vspace{1ex}
\caption{New values of $N_q(g)$ obtained in this paper.} 
\label{Tbl-exact}
\end{table}

Our improved bounds are due to the fact that some zeta 
functions that satisfy the numerical conditions of the 
explicit formul{\ae} method are forbidden by a combination of 
geometrical and numerical conditions.  Our approach is in the spirit of 
\cite{Lauter:CRAS,
Lauter:PAMS,
Lauter:conf,
LauterSerre:JAG,
LauterSerre:CM,
Savitt,
Serre:notes}
where lists of possible zeta functions were compiled and geometric arguments 
were applied for the purpose of improving the bounds. 

The main theorem that we use to improve the known upper bounds deals with 
a numerical invariant of pairs of abelian varieties.  Suppose $A_1$ and $A_2$ 
are abelian varieties over~$\F_q$.  Let $F$ and $V$ denote, respectively, 
the Frobenius and Verschiebung endomorphisms of $A_1\times A_2$. Given 
an element $\alpha$ of the subring $\Z[F,V]$ of $\End(A_1\times A_2)$, 
we let $g_1$ and $g_2$ be the minimal polynomials of $\alpha$ restricted 
to $A_1$ and $A_2$, respectively, and we define $r(\alpha)$ to be the 
resultant of $g_1$ and $g_2$. Define $s(A_1,A_2)$ to be the greatest 
common divisor of the set $\{r(\alpha) : \alpha \in \Z[F,V]\}$. 
Note that if $A_1$ and $A_2$ have an isogeny factor in common then 
$r(\alpha) = 0$ for every $\alpha$, so that $s(A_1,A_2) = \infty.$ 
On the other hand, if $A_1$ and $A_2$ share no common isogeny factors then
$r(F) \neq 0 $ by the Honda-Tate theorem~\cite{Tate}, so that 
$s(A_1,A_2) < \infty$. In other words, $\Hom(A_1,A_2) = \{0\}$ 
if and only if $s(A_1,A_2) < \infty$.
Also note that the value of $s(A_1,A_2)$ depends only on the
isogeny classes of $A_1$ and~$A_2$.

\begin{thm} 
\label{T-resultant} 
Let $A_1$ and $A_2$ be nonzero abelian varieties over $\F_q$. 
\begin{itemize} 
\item[(a)] If $s(A_1,A_2) = 1$ then there is no curve $C$ over $\F_q$ 
whose Jacobian is isogenous to $A_1\times A_2$. 
\item[(b)] Suppose $s(A_1,A_2) = 2$.  If $C$ is a curve over $\F_q$ 
whose Jacobian is isogenous to $A_1\times A_2$, then there is a 
degree-$2$ map from $C$ to another curve $D$ over $\F_q$ whose 
Jacobian is isogenous to either $A_1$ or~$A_2$. 
\end{itemize} 
\end{thm}

One can get upper and lower bounds on $s(A_1,A_2)$ by using the following
result.  (Recall that the {\em radical\/} of a nonzero integer is the 
product of its prime divisors.)

\begin{thm}
\label{T-radical}
Suppose $A_1$ and $A_2$ are abelian varieties over $\F_q$ with 
$s(A_1,A_2)\neq 0$.  Then $s(A_1,A_2)$ divides $r(F+V)$ and 
is divisible by the radical of~$r(F+V)$. 
\end{thm}

Theorem~\ref{T-radical} shows that Theorem~\ref{T-resultant}(a) 
is equivalent to a result of Serre~\cite{Lauter:PAMS, Serre:notes} 
that states that the Jacobian of a 
curve is never isogenous to a product $A_1\times A_2$ of
nonzero abelian varieties for which $r(F+V) = \pm 1$. 

It is not clear whether there are any similarly strong conclusions to be 
drawn from other values of $s(A_1,A_2)$.  However, if we make some 
assumptions about $A_1$ and $A_2$, we can prove that certain 
other values of $s(A_1,A_2)$ prohibit the existence of a curve with 
Jacobian isogenous to $A_1\times A_2$ --- see Propositions~\ref{P-general}
and~\ref{P-ec} and Corollaries~\ref{C-veryspecial} 
and~\ref{C-ec}.

Theorem~\ref{T-resultant}, combined with previously known results
and some straightforward facts about degree-$2$ maps of curves, 
allows us to greatly restrict the possible zeta functions of curves
having a large number $N$ of points.   For some values of $q$ and $g$
these restrictions are strong enough to allow us to immediately 
eliminate certain values of $N$ from consideration.
For other combinations of $q$, $g$, and $N$, 
we can quickly eliminate most possible zeta functions and
are left with a few special cases to consider. 
For some of these special cases we can use Theorem~\ref{T-resultant}(b)
to restrict the form of the curves in question to such an extent that
a computer search for curves with the desired number of points becomes
feasible.  For one such case, detailed in Section~\ref{S-descent},
we manage to avoid significant computer calculations by extending a 
Galois descent argument used in~\cite{Serre:notes}.

The {\em defect\/} of a genus-$g$ curve $C$ over $\F_q$ 
is the difference between the Weil-Serre upper bound 
for genus-$g$ curves over $\F_q$ and the number 
of rational points on~$C$; in other words, a curve $C$ has defect $k$ 
if it has exactly $(q + 1 + g[2\sqrt{q}]) - k$ rational points.
Theorem~\ref{T-resultant} allows us to prove some general
results about curves with small defect.
For example, we have the following theorem for square $q$.
\begin{thm} 
\label{T-square}
Suppose $q$ is a square.
\begin{itemize} 
\item[(a)] If $q \ne 4$ and $g>2$ then there are no defect-$2$ curves
           of genus $g$ over $\F_q$.
\item[(b)] If $q \ne 9$ and $g>3$ then there are no defect-$3$ 
           curves of genus $g$ over $\F_q$.
\item[(c)] If $g > (3q + 4m - 9) / m$, where $m = 2\sqrt{q}$,
           then there are no defect-$4$ curves of genus $g$ over $\F_q$.    
\item[(d)] If $q=2^{2e}$ with $e>2$, and if $g > 2^{e-1} + 2$, then there are
           no defect-$4$ curves of genus $g$ over $\F_q$.
\end{itemize} 
\end{thm}

For certain nonsquare $q$ the Weil-Serre bound can
be improved via a different method. 
Suppose $q$ is a prime power.
We define the {\em defect-$0$ dimension\/} of $q$ to be the
smallest positive integer $\delta$ for which there is a $\delta$-dimensional
abelian variety over $\F_q$ with characteristic polynomial
of Frobenius equal to $(x^2 + mx + q)^\delta$.
We say that $q$ is {\em exceptional\/} if its defect-$0$ dimension
is greater than~$1$.

\begin{thm}
\label{T-exceptional}
Suppose $q$ is a prime power and let $\delta$ be the
defect-$0$ dimension of~$q$. If $C$ is a curve of genus $g$ over $\F_q$, then
the defect of $C$ is at least $r/2$, where $r\in [0,\delta)$ is 
the remainder when $g$ is divided by~$\delta$.
\end{thm}

Theorem~\ref{T-exceptional} says something nontrivial about $q$ only 
if $q$ is exceptional, so we would like
to be able to find the exceptional prime powers.
In fact, there is an easy way to calculate the defect-$0$ dimension of
a power $q$ of a prime~$p$.  Let $\nu$ be an additive $p$-adic
valuation on $\Q$ and let $m = [2\sqrt{q}]$.

\begin{prop}
\label{P-defect-zero-dimension}
If $q$ is a square or if $q<4$ then the defect-$0$ dimension of $q$ is~$1$.
If $q>4$ is not a square, then the defect-$0$ dimension of $q$ is the
smallest positive integer $\delta$ such that $\delta \nu(m)/\nu(q)$ is an integer.
\end{prop}

We will prove Theorem~\ref{T-exceptional} and 
Proposition~\ref{P-defect-zero-dimension} in Section~\ref{S-exceptional}.  
The proofs will foreshadow the arguments that produce the entries in
Table~\ref{Tbl-char2} for $q = 128$.

It is easy to show that there are infinitely many $q$ of the form
$2^{2e+1}$ whose defect-$0$ dimension is~$2e+1$; we will provide a
proof of this fact in Section~\ref{S-exceptional}.
For such a $q$ we see that a curve of genus $g\le 2e$ must have defect at least~$g/2$.
The existence of these $q$ allows us to prove an interesting fact about
the function $N_q(g)$ defined by
$$N_q(g) = \max \{\#C(\F_q) : \text{$C$ is a genus-$g$ curve over $\F_q$}\}.$$

\begin{cor}
\label{C-exceptional}
There are infinitely many powers $q$ of $2$ such that for
every $g$ with $0 < g < \log_2 q$ we have
$(q + 1 + g [2\sqrt{q}]) - N_q(g) \ge g/2.$
\end{cor}

In particular, this implies that there is a sequence of pairs $(q,g)$ 
where $g$ is small with respect to $q$ and for which the Weil-Serre
bound becomes arbitrarily far from the true value of $N_q(g)$.
Zieve~\cite{Zieve} has already shown that there is a sequence 
of pairs $(q,g)$ where $g/q \to 1/2$ and for which
all previously-known bounds on $N_q(g)$ become 
arbitrarily far from the true value of $N_q(g)$.

Savitt~\cite{Savitt} recently showed, through extensive computer 
calculation, that there is no genus-$4$ curve over $\F_8$ having
exactly $27$ rational points. We prove this same result in 
Section~\ref{S-Hermitian} with an argument much less dependent
on the computer.  In~\cite{LauterSerre:JAG} it was
shown that there are only two possibilities for the
zeta function of such a curve. 
We can show that the first zeta function cannot occur by 
using Theorem~\ref{T-resultant}(b).  For the second zeta function, 
we introduce a new argument that generalizes the Hermitian form 
argument used in~\cite{LauterSerre:CM}.  We are able to eliminate 
the second zeta function by showing that if $A$ is an abelian 
variety whose characteristic polynomial of Frobenius is $f^2$, 
where 
$$f = x^4 - 9x^3 + 35x^2 - 72x + 64,$$
then every principal polarization on $A$ is decomposable.
To prove this, we show that there are no indecomposable 
unimodular Hermitian forms of rank~$2$ over the ring of integers
of the quartic number field defined by~$f$.

In Section~\ref{S-double} we prove Theorems~\ref{T-resultant}, 
\ref{T-radical}, and \ref{T-square},
and we provide a number of 
useful corollaries.  In Section~\ref{S-exceptional}
we prove Theorem~\ref{T-exceptional}, 
Proposition~\ref{P-defect-zero-dimension}, and
Corollary~\ref{C-exceptional}.
In Section \ref{S-improvements} we prove 
the results mentioned in Tables~\ref{Tbl-char2} 
and~\ref{Tbl-char3}, although we postpone 
the consideration of some cases to later sections. 
In Section~\ref{S-descent} we use a Galois descent argument to 
show that there is no genus-$5$ curve over $\F_8$ with $31$ points. 
In Section~\ref{S-exhaustion} we check by exhaustion that 
there is no genus-$4$ curve over $\F_{27}$ with $66$ points,
that there is no genus-$4$ curve over $\F_{32}$ with $75$ points,
and that there is no genus-$6$ curve over $\F_3$ with a certain
Weil polynomial; 
these calculations are feasible only because Theorem~\ref{T-resultant} 
allows us to considerably reduce the spaces we must search over. 
In Section~\ref{S-char3} we show how one can parametrize 
degree-$3$ covers of elliptic curves in characteristic~$3$.
Originally we had hoped to use this parameterization to show
show that there is no genus-$6$ curve over $\F_3$ with a certain 
Weil polynomial and that there is no genus-$4$ curve over $\F_{27}$
with $65$ points; however, our original arguments were flawed.
We sketch correct arguments in the second Appendix.
In Section~\ref{S-Hermitian} we use the Hermitian form argument 
mentioned above to prove Savitt's result that 
there is no genus-$4$ curve over $\F_8$ with $27$ points.

\subsubsection*{Notation} 
Throughout this paper a {\it curve} 
over $\F_q$ will mean a smooth, projective, absolutely irreducible curve. 
We will denote by $N_q(g)$ the largest $N$ such that there is a curve of 
genus $g$ over $\F_q$ with exactly $N$ rational points.
The {\em Weil polynomial\/} of an abelian variety over a finite field
is the characteristic polynomial of the Frobenius endomorphism of the
variety.  The {\em Weil polynomial\/} of a curve over a finite field
is the Weil polynomial of its Jacobian.  Note that if $f\in\Z[x]$ is 
the Weil polynomial of a genus-$g$ curve $C$ over $\F_q$, then the
numerator of the zeta function of $C$ is equal to~$x^{2g} f(1/x).$
If $f$ is the Weil polynomial of a curve or an abelian variety, say 
with $\deg f = 2g$, then there is a degree-$g$
polynomial $h\in\Z[x]$, all of whose roots are
real, such that $f(x) = x^g h(x + q/x).$
We will refer to $h$ as the {\em real Weil polynomial\/} of the
curve or the abelian variety.

\subsubsection*{Acknowledgments}
We thank Jean-Pierre Serre for his helpful comments.
In the course of doing the work described in this paper
we used the computer algebra systems Pari/GP and Magma~\cite{magma}.
Several of our Magma programs are available on the web:  start at

\smallskip
\noindent
{\tt http://www.alumni.caltech.edu/\~{}however/biblio.html}
\smallskip

\noindent
and follow the links related to this paper.  We have also placed
a copy of the 18 January 2002 version of the van der Geer-van der Vlugt
tables on this site.

\section{Proofs of Theorems~{\rm\ref{T-resultant}}, 
         {\rm\ref{T-radical}}, and {\rm\ref{T-square}}}
\label{S-double}

In this section we will prove Theorems~\ref{T-resultant}, 
\ref{T-radical}, and~\ref{T-square},
as well as some useful corollaries.  We begin with a lemma.

\begin{lem}
\label{L-resultant}
Suppose $B$ is an abelian variety over $\F_q$ isogenous to 
a product $A_1\times A_2$, where $s(A_1,A_2) < \infty$.
Then there exist abelian varieties
$A_1'$ and $A_2'$, isogenous to $A_1$ and $A_2$, respectively,
and an exact sequence 
$$0\to \Delta' \to A_1'\times A_2' \to B \to 0$$ 
such that the projection maps $A_1'\times A_2' \to A_1'$ 
and $A_1'\times A_2' \to A_2'$ give monomorphisms
from $\Delta'$ to $A_1'[s]$ and to $A_2'[s]$,
where $s = s(A_1,A_2)$.

Suppose in addition that $B$ has a principal polarization~$\mu$.
Then the pullback of $\mu$ to $A_1' \times A_2'$ is a
product polarization $\lambda_1 \times \lambda_2$, and 
the projection maps $A_1'\times A_2' \to A_1'$ 
and $A_1'\times A_2' \to A_2'$ give isomorphisms of $\Delta'$
with $\ker \lambda_1$ and $\ker \lambda_2$.
In particular, $\Delta'$ is isomorphic to its own Cartier dual.
\end{lem}

\begin{proof}
Let $\varphi$ be an arbitrary isogeny from $A_1\times A_2$ to $B$ and let
$\Delta$ be the kernel of $\varphi$.  Let $G_1$ and $G_2$ be the largest
closed subgroup-schemes of $A_1$ and $A_2$ such that $G_1\times G_2$ is a 
closed subgroup-scheme of $\Delta$, let $\Delta' = \Delta / (G_1\times G_2)$, 
and let $A_1' = A_1/G_1$ and $A_2' = A_2 / G_2$.
Then we have an exact sequence
$$0\to \Delta' \to A_1'\times A_2' \to B \to 0$$
such that the projection maps give monomorphisms of $\Delta'$ to $A_1'$ 
and $A_2'$.  We will show that in fact the projection maps take $\Delta'$
to $A_1'[s]$ and $A_2'[s]$.

Let $\alpha$ be an arbitrary endomorphism of 
$A_1' \times A_2'$ that lies in $\Z[F,V]$.
(Here we use the fact that $\Z[F,V]$ is contained in the
endomorphism ring of every abelian variety isogenous to $A_1\times A_2$.)
For $i=1,2$ let $g_i$ be the minimal polynomial of $\alpha$ acting on $A_i'$.
Then $g_1(\alpha)$ and $g_2(\alpha)$ both act as $0$ on $\Delta'$, because 
$\Delta'$ can be viewed as a subgroup-scheme of both $A_1'$ and $A_2'$. 
But then $r(\alpha)$ must act as $0$ on $\Delta'$ as well. Since this
is true for every $\alpha$, we see that $s(A_1,A_2)$ must kill~$\Delta'$;
this shows that the projection maps embed $\Delta'$ into $A_1'[s]$ 
and~$A_2'[s]$. 

Now suppose $B$ has a principal polarization~$\mu$.
Let $\lambda$ be the pullback of $\mu$ to $A_1'\times A_2'$.
Since $\Hom(A_1', A_2')$ and $\Hom(A_2',A_1')$ are both trivial,
$\lambda$ must be a product polarization $\lambda_1\times\lambda_2$.
The degree of $\lambda$ is equal to the degree of $\mu$ (which is~$1$)
times the square of the degree of the isogeny $A_1'\times A_2'\to B$,
so we have
$$(\#\Delta')^2 = \#\ker\lambda = (\#\ker\lambda_1)(\#\ker\lambda_2),$$
where we use $\#$ to denote the rank of a finite group-scheme.
Since the projection maps give monomorphisms from $\Delta'$
to $A_1'$ and $A_2'$, we see that $\#\Delta' \le \#\ker\lambda_i$
for $i=1$ and $i=2$.  This means that we must have 
$\#\Delta' = \#\ker\lambda_i$ for each $i$, and it follows that
$\Delta' \cong \ker\lambda_i$ for each $i$.  Since kernels of 
polarizations are isomorphic to their own duals, we obtain the 
final statement of the lemma.
\end{proof}

\begin{proof}[Proof of Theorem~{\rm\ref{T-resultant}}]
Suppose $s(A_1,A_2) = 1$.  Then Lemma~\ref{L-resultant} shows that every
abelian variety isogenous to $A_1\times A_2$ is a product $A_1'\times A_2'$.
Since $s(A_1',A_2') = s(A_1,A_2) <\infty$ we see that 
$\Hom(A_1',A_2') = \{0\}$, so every polarization on $A_1'\times A_2'$ 
is a product polarization.  In particular, we see that every principal
polarization of an abelian variety isogenous to $A_1\times A_2$ is 
decomposable, so there can be no Jacobians isogenous to $A_1\times A_2$.
This is the first statement of the theorem.

Now suppose that $s(A_1,A_2) = 2.$  Apply Lemma~\ref{L-resultant} and
replace $A_1$ and $A_2$ with the resulting $A_1'$ and $A_2'$, so that 
we have an exact sequence 
$$0\to \Delta \to A_1\times A_2 \to \Jac C \to 0$$
where $\Delta$ can be viewed as a subscheme of $A_1[2]$ and $A_2[2]$.

Let $\mu$ be the canonical polarization on $\Jac C$ and let $\lambda$ 
be the polarization on $A_1\times A_2$ obtained by pulling back $\mu$ 
via~$\varphi$.  Lemma~\ref{L-resultant} shows that $\lambda$
is the product of a polarization $\lambda_1$ on $A_1$ and a 
polarization $\lambda_2$ of $A_2$.
Let $(1,-1)$ denote the involution of $A_1\times A_2$ that acts 
as $1$ on $A_1$ and as $-1$ on $A_2$.  Clearly $(1,-1)$ respects the 
polarization $\lambda$, because $1$ respects $\lambda_1$ and
$-1$ respects $\lambda_2$. 
Furthermore, $(1,-1)$ acts as the identity on $\Delta$, so it gives rise 
to an involution $\beta$ on $\Jac C$ that respects the
polarization~$\mu$. By Torelli's theorem, there exists 
an involution $\alpha$ of $C$ such that either $\beta = \alpha^*$ or 
$\beta = -\alpha^*$.  

Let $D$ be the quotient of $C$ by the involution $\alpha$, so that 
there is a degree-$2$ map $\psi$ from $C$ to $D$ with 
$\psi = \psi\circ\alpha$. Then the morphism $\psi^*:\Jac D\to \Jac C$ 
gives an isogeny from $\Jac D$ to the connected component of the 
subvariety of $\Jac C$ on which $\beta$ acts as the identity.  
This subvariety is isogenous to $A_1$ if $\beta = \alpha^*$ and 
to $A_2$ if $\beta = -\alpha^*$. 
\end{proof}

We will use Theorem~\ref{T-resultant}(b) in conjunction with
some obvious facts about degree-$2$ covers of curves, which we
state here for convenience.

\begin{lem}
\label{L-obvious}
Suppose $C$ and $D$ are curves over $\F_q$ of genus $g_C$ and $g_D$,
respectively, and suppose there is a degree-$2$ map $C\to D$.
For every integer $d>0$ let $a_d$ denote the number of degree-$d$
places on $C$ and let $b_d$ denote the number of degree-$d$
places on $D$.
\begin{itemize}
\item[(a)] For every odd $d$ we have $a_d\le 2 b_d$.
\item[(b)] We have $g_C \ge 2g_D - 1$, with equality
           if and only if $C\to D$ is unramified.
\item[(c)] Let $d_1 < \cdots < d_n$ be odd positive integers
           such that $a_{d_i}$ is odd for every~$i$, and let
           $r = d_1 + \cdots + d_n$.  Then
           $g_C \ge 2g_D - 1 + r/2,$ and equality holds
           if and only if $C\to D$ is ramified at 
           exactly $n$ places $\frakp_1,\ldots,\frakp_n$, where each
           $\frakp_i$ has degree $d_i$ and where the
           ramification at each $\frakp_i$ is tame.
\end{itemize}
\end{lem}

\begin{proof}
Suppose $d$ is odd. 
Every degree-$d$ place of $D$ has at most $2$ degree-$d$ places of $C$ lying 
over it, and every degree-$d$ place of $C$ lies over a degree-$d$ place of $D$.
Statement (a) follows immediately.

Statement (b) is the special case $n=0$ of statement (c).

Let $\iota$ be the involution of $C$ corresponding to the cover $C\to D$.
Suppose $d$ is an odd number such that $a_d$ is odd.  Then there is
a degree-$d$ place $\frakp$ of $C$ that is taken to itself by $\iota$.
Since $\frakp$ consists of an odd number of geometric points of~$C$,
there must be a geometric point $P$ in $\frakp$ that is fixed by~$\iota$,
and since all of the geometric points in $\frakp$ are conjugate to
each other, {\em all\/} of the points in $\frakp$ must be fixed by~$\iota$.
These $d$ points must be ramification points of the cover $C\to D$.
Thus, the hypothesis of statement (c) implies that there are at least $r$
ramification points in the cover~$C\to D$.  The conclusion of the
statement then follows by applying the Riemann-Hurwitz formula to
the cover~$C\to D$.
\end{proof}

Suppose $A$ is a $g$-dimensional abelian variety over~$\F_q$
and let $$\{\alpha_1,\ldots,\alpha_g, \alphabar_1,\ldots,\alphabar_g\}$$
be the multiset of complex roots of the Weil polynomial of~$A$.
For each $i$ let $x_i = -(\alpha_i + \alphabar_i)$. We will say that $A$ is 
{\em of type $[x_1,\ldots, x_g]$}. If $A$ is the Jacobian of a curve $C$ 
we will also say that $C$ and its zeta function are 
{\em of type $[x_1,\ldots, x_g]$}. Note that the zeta function of a curve 
of type $[x_1,\ldots, x_g]$ is given by 
$$\frac{n(t)}{(1-t)(1-qt)},$$ 
where 
$$n(t) = \prod_{i=1}^g (1 + x_i t + qt^2).$$
Also, if $F$ is the Frobenius morphism of $A$ and if $V = q/F$ is
the Verschiebung, then the characteristic polynomial of $F+V$ is
equal to $h^2(t)$, where
$$h(t) = \prod_{i=1}^g (t + x_i)$$
is the real Weil polynomial of~$A$.

\begin{cor}
\label{C-type2}
There are no genus-$g$ curves of type $[m,\ldots,m,m-2]$ over $\F_q$ 
if  $g > (q-1+2m)/m$ and $g>3$, where $m = [2\sqrt{q}]$. 
\end{cor}

\begin{proof}
We will prove the contrapositive statement.
Suppose $C$ is a curve of genus $g$ over $\F_q$  with zeta function 
$[m,\ldots,m,m-2]$. Then $\Jac C$ is isogenous to a product $A\times E$ 
of abelian varieties, where $E$ is an elliptic curve over $\F_q$ of 
type $[m-2]$ and where $A$ is a $(g-1)$-dimensional abelian variety 
over $\F_q$ of type $[m,\ldots,m]$. We see that $r(F+V) = 2$, so 
$s(A,E) = 2$. According to Theorem \ref{T-resultant}, the curve $C$ 
is a degree-$2$ cover of a curve $D$ whose Jacobian is isogenous to 
either $A$ or $E$. If $\Jac D \sim A$ then $D$ has genus $g-1$, and 
Lemma~\ref{L-obvious}(b) shows that~$g \le 3$. 
If $\Jac D\sim E$ then $D$ is an elliptic 
curve with $q + m - 1$ points, and applying Lemma~\ref{L-obvious}(a)
with $d=1$ shows that 
$$q + gm - 1\le 2q + 2m - 2,$$ which gives $g \le (q-1+2m)/m$. 
\end{proof}

Recall that the {\em defect\/} of a genus-$g$ curve $C$ over $\F_q$ 
is the difference between the Weil-Serre upper bound and the number 
of rational points on~$C$.

\begin{cor}
\label{C-defect2}
There are no defect-$2$ curves of genus $g$ over $\F_q$ if
$g > (q-1+4m)/m$ and $g>5$, where $m = [2\sqrt{q}]$.
\end{cor}

\begin{proof}
If $C$ has defect $2$, then its zeta function must be of one of the 
seven types listed in~\cite{LauterSerre:JAG}.  For $g \ge 5$, all but 
two of these types are eliminated by Theorem~\ref{T-resultant}(a). 
The two remaining types are $[m,\ldots,m,m-2]$ and
$[m,\ldots,m,m+\sqrt{3}-1,m-\sqrt{3}-1]$. 
Since we are assuming that $g > (q-1+4m)/m$, Corollary~\ref{C-type2}
eliminates  the former possibility, so $C$ must have the
latter type. In this
case $\Jac C$ is isogenous to the product of a $(g-2)$-dimensional 
abelian variety $A_1$ of type $[m,\ldots,m]$ and a $2$-dimensional
abelian variety $A_2$ of type $[m+\sqrt{3}-1, m-\sqrt{3}-1]$. 
Applying Theorem~\ref{T-resultant}(b), we find that $C$ is a 
double cover of a curve $D$ that is either of type $[m,\ldots,m]$ 
or of type $[m+\sqrt{3}-1,m-\sqrt{3}-1]$. In the first case $D$ 
would have genus $g-2$, and Lemma~\ref{L-obvious}(b)
shows that then~$g\le 5$.  In the second case, 
Lemma~\ref{L-obvious}(a) with $d=1$ gives us
$$q + gm - 1 \le 2 (q + 2m - 1),$$
which leads to $g \le (q - 1 + 4m)/m$.
\end{proof}

We have mentioned that we do not know any strong conclusions one 
can draw in general when $s(A_1,A_2) > 2$.  However, with a little 
more information about $A_1$ and $A_2$ we can indeed say something.

\begin{prop}
\label{P-general}
Let $A_1'$ and $A_2'$ be abelian varieties over $\F_q$ and let 
$s = s(A_1',A_2')$. Suppose that for every $A_1$ isogenous to $A_1'$
and every $A_2$ isogenous to $A_2'$, the only self-dual finite 
group-scheme that can be embedded in both $A_1[s]$ and $A_2[s]$ 
as the kernel of a polarization
is the trivial group-scheme.  Then there is no curve over $\F_q$ 
with Jacobian isogenous to $A_1'\times A_2'$. 
\end{prop}

\begin{proof}
Suppose there were such a curve.
Then Lemma~\ref{L-resultant} shows that we can find a group-scheme 
$\Delta$ that fits in an exact sequence 
$$0\to \Delta \to A_1\times A_2 \to \Jac C \to 0$$
and that can be embedded in both $A_1[s]$ and $A_2[s]$. 
Furthermore, since $\Jac C$ has a principal polarization, 
for each $i=1,2$ we have that $\Delta$ is isomorphic to the 
kernel of the polarization of $A_i$ obtained by pulling back the
principal polarization of~$\Jac C$.
The hypotheses of the proposition show that $\Delta$ must be the 
trivial group-scheme, so $\Jac C$ is isomorphic to the product of
two abelian varieties that share no isogeny factor.  
As we have seen, this is a contradiction. 
\end{proof}

The next corollary describes a situation in which the hypotheses 
of Proposition~\ref{P-general} are met.

\begin{cor}
\label{C-veryspecial}
Suppose $q$ is a square prime power and $n$ is a squarefree integer
coprime to~$q$. 
Let~$m = 2\sqrt{q}$.  Then there is no curve over $\F_q$ 
of type $[m,\ldots,m,m-n]$. 
\end{cor}

\begin{proof}
Let $A_1$ be any abelian variety over $\F_q$
isogenous to the product of $(g-1)$ copies of a supersingular 
elliptic curve over $\F_q$ with Weil polynomial 
$x^2 + mx + q = (x + \sqrt{q})^2$, and let $A_2$ be any ordinary elliptic curve 
over $\F_q$ with Weil polynomial $x^2 + (m-n)x + q$. 
Clearly $s(A_1,A_2) = n$.  Suppose $\Delta$ is a nontrivial
self-dual finite group scheme that embeds in both $A_1[n]$
and $A_2[n]$ as the kernel of a polarization.
Since the only polarizations on $A_2$ are the
multiplication-by-$\ell$ maps for positive integers~$\ell$,
we must have $\Delta\cong A_2[\ell]$ for some divisor $\ell>1$ of $n$.
Since $\Delta$ embeds in $A_1[\ell]$
as well, and since the Frobenius $F$ on $A_1$ satisfies $F + \sqrt{q} = 0$, 
we know that Frobenius must act as the integer
$-\sqrt{q}$ on $\Delta$, and hence on $A_2[\ell]$.  This means
that $F + \sqrt{q} = 0$ on $A_2[\ell]$, which means that
$(F + \sqrt{q})/\ell$ is an endomorphism of $A_2$.  But
from the characteristic polynomial of $F$ on $A_2$ we can
calculate that the characteristic polynomial of $(F + \sqrt{q})/\ell$
on $A_2$ is 
$$x^2 - (n/\ell) x + (n\sqrt{q}/\ell^2),$$which is not integral.
This contradiction shows that no nontrivial self-dual finite group-scheme 
can be embedded in both $A_1[n]$ and $A_2[n]$.
By Proposition~\ref{P-general}, there is no curve over $\F_q$
of type $[m,\ldots,m,m-n]$. 
\end{proof}

There is another situation in which we can draw conclusions from
values of $s(A_1,A_2)$ greater than $2$.

\begin{prop}
\label{P-ec}
Suppose $C$ is a curve over $\F_q$ whose Jacobian is isogenous
to the product $A\times E$ of an abelian variety $A$ with an
elliptic curve $E$, and suppose that $s(A,E) < \infty$.
Then there is an elliptic curve $E'$ isogenous to $E$ for which
there is map from $C$ to $E'$ of degree dividing $s(A,E)$,
and we have $\#C(\F_q) \le s(A,E) \cdot \#E(\F_q)$.
\end{prop}

\begin{proof}
By applying Lemma~\ref{L-resultant} we find that there are
abelian varieties $A'$ and $E'$, isogenous to $A$ and $E$, 
respectively, and an exact sequence
$$0\to \Delta' \to A'\times E' \to \Jac C\to 0$$ 
such that the projection maps $A'\times E' \to A'$ 
and $A'\times E' \to E'$ give monomorphisms
from $\Delta'$ to $A'[s]$ and $E'[s]$,
where $s = s(A,E)$.  This implies that the
composition $E'\to A'\times E' \to \Jac C$
is a monomorphism.   Let $\lambda$ and $\mu$ be the canonical principal
polarizations on $E'$ and $\Jac C$, respectively.
Then the pullback of $\mu$ to $E'$ is $n\lambda$ for some integer~$n>0$,
and Lemma~\ref{L-resultant} says that the kernel of the pullback
is isomorphic to $\Delta'$.
In particular we see that $n$ must divide~$s$.
Thus we have a diagram
\begin{equation*}
\begin{matrix}
     E'     & \mapright{n} &     E'    & \mapright{\lambda} & \hat{E'}    \\
\Bdownarrow &              & \Buparrow &                    &\Buparrow    \\
  \Jac C    & \mapright{1} &   \Jac C  & \mapright{\mu}     & \hat{\Jac C}
\end{matrix}
\end{equation*}
where the vertical arrow on the right is the dual morphism of
the vertical monomorphism $E'\to\Jac C$ on the left.  From this 
we see that the composition $C\to\Jac C\to E'$ is a map of degree~$n$.
In particular,
$$\#C(\F_q)\le n \cdot\#E'(\F_q) = n\cdot\#E(\F_q) \le s \cdot\#E(\F_q).$$
\end{proof}

\begin{cor}
\label{C-ec}
Let $q$ be a prime power and let $n$ be a positive integer.
Let~$m = \lfloor 2\sqrt{q}\rfloor$.  If $C$ is a curve over $\F_q$
of type $[m,\ldots,m,m-n]$, then the genus $g$ of $C$ satisfies
$$g \le \frac{(n-1)q - (n-1)^2 + nm}{m}.$$
\end{cor}

\begin{proof}
According to Proposition~\ref{P-ec}, there is an elliptic curve $E'$
of defect $n$ and a divisor $d$ of $n$ for which
there is a degree-$d$ map from $C$ to $E'$, and the number
of points on~$C$, which is $q+1+gm - n$, is at most $n$ times 
the number of points on~$E'$, which is $q + 1 + m - n$.
It follows from this inequality that 
$$g \le \frac{(n-1)q - (n-1)^2 + nm}{m}.$$
\end{proof}

Now we turn to the proof of Theorem~\ref{T-radical}.  Our proof relies 
on two facts: First, a prime $p$ 
divides the resultant of two monic polynomials in $\Z[x]$ if and only 
if the reductions of the polynomials modulo $p$ have a common root in 
the algebraic closure of~$\F_p$, and second, the ring $\Z[F,V]$ contains 
no nilpotent elements, so that if $\beta\in\Z[F,V]$ has minimal 
polynomial $m\in\Z[x]$, then the subring $Z[\beta]$ of $\Z[F,V]$ is 
isomorphic to $\Z[x] / (m)$.

\begin{proof}[Proof of Theorem~{\rm\ref{T-radical}}]
Clearly $s(A_1,A_2)$ divides $r(F+V)$, because $s(A_1,A_2)$ is defined 
to be the greatest common divisor of a set of numbers that includes
$r(F+V)$.
What we must now prove is that if a prime $p$ divides $r(F+V)$ then 
it also divides $s(A_1,A_2)$.  To do this, we must show that for 
every $\alpha$ in $\Z[F,V]$ the prime $p$ divides~$r(\alpha)$.

Consider an $\alpha$ in $\Z[F,V]$, say $\alpha = u(F,V)$ for some
polynomial $u\in\Z[x,y]$.   
Let $F_1$ and $V_1$ (resp.~$F_2$ and $V_2$) 
be the Frobenius and Verschiebung endomorphisms of $A_1$ (resp.~$A_2$). 
The fact that $p$ divides $r(F+V)$ shows that there are homomorphisms 
$\psi_1\colon\Z[F_1+V_1] \to \Fpbar$ and $\psi_2\colon\Z[F_2+V_2] \to \Fpbar$ 
such that $\psi_1(F_1+V_1) = \psi_2(F_2+V_2)$.  Let $\tau = \psi_1(F_1+V_1)$ 
and let $\sigma$ be an element of $\Fpbar$ such that 
$\sigma^2 - \tau\sigma + q = 0$.  Then the homomorphisms $\psi_1$ and 
$\psi_2$ can be extended to give homomorphisms 
$\psihat_1\colon\Z[F_1,V_1]\to\Fpbar$ and 
$\psihat_2\colon\Z[F_2,V_2]\to\Fpbar$
such that $\psihat_1(F_1) = \sigma = \psihat_2(F_2)$
and $\psihat_1(V_1) = \tau - \sigma = \psihat_2(V_2)$.
But then $\psihat_1(u(F_1,V_1)) = \psihat_2(u(F_2,V_2)),$
so the minimal polynomials of $u(F_1,V_1)$ and $u(F_2,V_2)$ have 
a common root in~$\Fpbar$.  It follows that $r(\alpha)$ is divisible 
by~$p$.
\end{proof}

Using Theorem~\ref{T-resultant} and the corollaries established 
so far, we can now prove Theorem~\ref{T-square}.  

\begin{proof}[Proof of Theorem~{\rm\ref{T-square}}]
Suppose $q$ is a square prime power, say $q=p^{2e}$ for
a prime~$p$.  Of the types of defect-$2$ zeta-functions listed 
in~\cite{LauterSerre:JAG}, only two are possible when $q$ is
a square: namely $[m,\ldots,m,m-1,m-1]$ and $[m,\ldots,m,m-2]$.
The first of these is impossible when $g>2$ by Theorem~\ref{T-resultant}(a),
so the only possible defect-$2$ zeta function for $g>2$ is $[m,\ldots,m,m-2]$.
To prove part (a), first assume that $p \ne 2$.  Then it follows 
from Corollary~\ref{C-veryspecial} that this zeta function is not 
possible.  
If $p=2$, then $m-2$ is not the trace of an elliptic curve when~$q\ne 4$,
so $[m,\ldots,m,m-2]$ is not possible in that case either.

The proof of part (b) is essentially the same. From~\cite{Savitt}
we see that the only possible zeta function for a defect $3$ curve when 
$q$ is a square and $g>3$ is $[m,\ldots,m,m-3]$.  
(If $g=3$, then $[m-1,m-1,m-1]$ may be possible; see~\cite{LauterSerre:CM}.)  
If $p \ne 3$ then it again follows from Corollary~\ref{C-veryspecial}
that this zeta function is not possible.  
If $p=3$, then $m-3$ is not the trace of an elliptic curve when 
$q\ne 9$, so $[m,\ldots,m,m-3]$ is not possible in that case either.

Now we prove parts (c) and~(d).  
Using the methods of~\cite{Serre:notes} (see also~\cite[\S2]{LauterSerre:JAG})
and the tables from~\cite{Smyth},
we see that for square $q$
there are exactly eight possible types for a curve of genus $g$
and defect $4$ over $\F_q$.
For each type, we list in Table~\ref{Tbl-defect4} 
the associated real Weil polynomial $h$ evaluated at $x - m$ 
(where $m = 2\sqrt{q}$)
and the associated Weil polynomial~$f$.
\begin{table}
\renewcommand\arraystretch{1.25}
\begin{center}
\begin{tabular}{|l|}
\hline  
type $=[m,\ldots,m,m-4]$                                     \\
$h(x-m) = x^{g-1} \cdot (x-4) $                              \\ 
$f(x) = (x + \sqrt{q})^{2g-2} \cdot (x^2 + (m-4)x + q)$      \\ \hline
type $=[m,\ldots,m,m-2,m-2]$                                 \\
$h(x-m) = x^{g-2}\cdot (x-2)^2 $                             \\ 
$f(x) = (x + \sqrt{q})^{2g-4} \cdot(x^2 + (m-2)x + q)^2$     \\ \hline
type $=[m,\ldots,m,m-1,m-1,m-2]$                             \\
$h(x-m) = x^{g-3} \cdot (x-1)^2\cdot (x-2)$                  \\
$f(x) = (x + \sqrt{q})^{2g-6} \cdot (x^2 + (m-1)x + q)^2 
                           \cdot (x^2 + (m-2)x + q)$         \\ \hline
type $=[m,\ldots,m,m-(2-\sqrt{2}),m-(2+\sqrt{2})]$           \\
$h(x-m) = x^{g-2}\cdot (x^2 - 4x + 2) $                      \\
$f(x) = (x + \sqrt{q})^{2g-4} \cdot 
     (x^4 + (2m-4)x^3 + (6q-4m+2)x^2 + (2m-4)qx + q^2)$      \\ \hline
type $=[m,\ldots,m,m-1,m-1,m-1,m-1]$                         \\
$h(x-m) = x^{g-4}\cdot (x-1)^4$                              \\
$f(x) = (x + \sqrt{q})^{2g-8} \cdot(x^2 + (m-1)x + q)^4$     \\ \hline
type $=[m,\ldots,m,m-1,m-(3-\sqrt{5})/2,m-(3+\sqrt{5})/2]$   \\
$h(x-m) = x^{g-3} \cdot(x-1) \cdot(x^2-3x+1)$                \\
$f(x) = (x + \sqrt{q})^{2g-6}\cdot (x^2 + (m-1)x + q) \cdot$ \\
$\qquad\qquad (x^4 + (2m-3)x^3 + (6q-3m+1)x^2 +(2m-3)qx + q^2)$    \\ \hline
type $=[m,\ldots,m,m-(2-\sqrt{3}),m-(2+\sqrt{3})]$           \\
$h(x-m) = x^{g-2}\cdot (x^2 - 4x + 1)$                       \\
$f(x) = (x + \sqrt{q})^{2g-4}\cdot 
     (x^4 + (2m-4)x^3 + (6q-4m+1)x^2 + (2m-4)qx + q^2)$      \\ \hline
type $=[m,\ldots,m,m-1,m-3]$                                 \\
$h(x-m) = x^{g-2}\cdot(x-1)\cdot(x-3)$                       \\
$f(x) = (x + \sqrt{q})^{2g-4} \cdot(x^2 + (m-1)x + q)
                           \cdot(x^2 + (m-3)x + q)$          \\ \hline
\end{tabular}
\end{center}
\vspace{1ex}
\caption{The possible types of defect-$4$ curves over square fields,
together with the associated real Weil polynomial $h$ evaluated at $x-m$
and the associated Weil polynomial~$f$.}
\label{Tbl-defect4}
\end{table}

Suppose that $g > (3q + 4m - 9)/m$, where $m = 2\sqrt{q}$.  
Then Corollary~\ref{C-ec} shows that the first entry in the table 
cannot occur.  Also, since the inequality we are assuming implies
that $g > (q  + 4m - 3)/m$ as well, Theorem~\ref{T-resultant}(b)
and Lemma~\ref{L-obvious}(a)
can be used to show that the second and fourth entries cannot occur.  
We also see that $g > (q + 2m + 3)/m$, so Theorem~\ref{T-resultant}(b)
and Lemma~\ref{L-obvious}(a) show that the eighth entry cannot occur.  
Finally, Theorem~\ref{T-resultant}(a)
shows that the remaining entries cannot occur when $g > 4$.
This proves part~(c).

Finally, suppose $q=2^{2e}$ with $e>2$, and suppose $g > 2^{e-1} + 2 > 5$.
Then the first four entries listed in Table~\ref{Tbl-defect4}
cannot occur because the Honda-Tate
theorem~\cite{Tate} shows the final factor of each of the
putative Weil polynomials is not in fact a Weil polynomial.
(A simple way to check this is to use~\cite[Lem.~3.1.2]{DiPippoHowe}.)
The next three entries cannot occur when $g > 4$ because
of Theorem~\ref{T-resultant}(a), as we have noted already.
That leaves us with the final entry.  Once again 
Theorem~\ref{T-resultant}(b) and Lemma~\ref{L-obvious}(a) 
can be used to eliminate this possibility,
because we have $g> 2^{e-1} + 2$.
\end{proof}

\section{Exceptional prime powers}
\label{S-exceptional}

In this section we will prove Theorem~\ref{T-exceptional}, 
Proposition~\ref{P-defect-zero-dimension}, and 
Corollary~\ref{C-exceptional}.  Before we begin,
let us define the {\em trace\/} of a monic degree-$n$
polynomial in $\Q[x]$ to be $-1$ times the coefficient
of $x^{n-1}$, and the {\em deficiency\/} of such a polynomial
to be its trace minus its degree.

\begin{proof}[Proof of Theorem~{\rm\ref{T-exceptional}}]
Let $C$ be a curve of genus $g$ over $\F_q$ and let $h\in\Z[x]$ be its
real Weil polynomial.  We know that all of the roots of $h$ are
real numbers in the interval $[-2\sqrt{q},2\sqrt{q}]$, and the
number of points on $C$ is equal to $q + 1 - t$, where $t$ is the
trace of~$h$.
Write $h = (x + m)^e h_2$, where $h_2$ has no factors
of $(x+m)$. 
The factor $(x+m)^e$ of $h$ corresponds to the largest
isogeny factor of $\Jac C$ on which Frobenius acts as~$-m$,
and up to isogeny this factor must be a power of the
smallest abelian variety over $\F_q$ whose real Weil polynomial
is a power of $x+m$.  Thus, the exponent $e$ is a multiple of the 
defect-$0$ dimension $\delta$ of $\F_q$, and we see that 
the degree $g_2$ of $h_2$ is congruent to $g$ modulo~$\delta$.

Define $H\in\Z[x]$ by $H(x) = h(x-m-1)$, so that $H = (x-1)^e H_2$
for a polynomial $H_2$ of degree $g_2$ that has no factors of $x-1$.
All of the roots of $H$ are positive real numbers, and the number of points
on $C$ is $q + 1 - T + gm + g$, where $T$ is the trace of~$H$.
The trace $T_2$ of $H_2$ is equal to $T - e$, and the degree 
of $h_2$ is equal to $g - e$, so we have
$$\#C(\F_q) = q + 1 + gm - T_2 + g_2.$$
In other words, the defect of the curve $C$ is $T_2 - g_2$, which is the
deficiency of the polynomial~$H_2$.

Now, all of the roots of $H_2$ are positive real numbers, and $H_2$ has
no factors of $x - 1$, so a result of Siegel~\cite{Siegel}
says that the trace of $H_2$ is at least $3/2$ times its degree.
It follows that the deficiency of $H_2$ is at least half its degree.
Thus, the defect of $C$ is at least $g_2/2$.  We already noted that 
$g_2$ is congruent to $g$ modulo~$\delta$,
so the defect of $C$ is at least $r/2$, where $r \in [0,\delta)$ is the
remainder obtained when dividing $g$ by~$\delta$.
\end{proof}

\begin{proof}[Proof of Proposition~{\rm\ref{P-defect-zero-dimension}}]
If $f$ is a monic irreducible polynomial in $\Z[x]$ whose 
roots in the complex plane all have magnitude $\sqrt{q}$, then
there is an exponent $e$ such that $f^e$ is the Weil polynomial
of a simple abelian variety over~$\F_q$.
The Honda-Tate theorem~\cite{Tate} includes a recipe for
calculating this exponent.  Proposition~\ref{P-defect-zero-dimension}
is obtained by applying this recipe to either the polynomial 
$x^2 + mx + q$ (if $q$ is not a square) or the polynomial $x + \sqrt{q}$
(if $q$ is a square).  We leave the details to the reader.
\end{proof}

The essence of the following argument appears in~\cite{Serre:notes}.

\begin{proof}[Proof of Corollary~{\rm\ref{C-exceptional}}]
Consider the expression for $\sqrt{2}$ in base~$2$:
$$\sqrt{2} = b_0 + \frac{b_1}{2} + \frac{b_2}{2^2} + \frac{b_3}{2^3} + \cdots$$
where each $b_i$ is either $0$ or $1$.
Suppose $e>0$ is an integer such that
$b_e = 1$ and $b_{e+1} = 0$.  
Let $q = 2^{2e+1}$.
Then the base-$2$ expression
for $2\sqrt{q}$ is
$$2\sqrt{q} = b_0 2^{e+1} + b_1 2^e + \cdots + b_e 2 + b_{e+1} 
                          + \frac{b_{e+2}}{2} +  \ldots,$$
so the base-$2$ expression for $m = [2\sqrt{q}]$ is
$$m  = b_0 2^{e+1} + b_1 2^e + \cdots +  b_e 2 + b_{e+1}.$$
Clearly $m$ is even but not a multiple of~$4$, 
so if $\nu$ is the usual additive $2$-adic valuation of $\Q$ 
we have $\nu(m) = 1$ and $\nu(q) = 2e+1$.  It follows 
from Proposition~\ref{P-defect-zero-dimension} that
the defect-$0$ dimension of $q$ is equal to $2e+1$.  If we take
$g \le 2e$, then Theorem~\ref{T-exceptional} shows that the
defect of a genus-$g$ curve over $\F_q$ is at least~$g/2$.  

Thus, to prove Corollary~\ref{C-exceptional} we need only show that
there are infinitely many $e$ with $b_e = 1$ and $b_{e+1} = 0$.
But there are only two ways in which there could not be
infinitely many such $e$: either $b_i = 0$ for all sufficiently
large $i$, or $b_i = 1$ for all sufficiently large $i$.  Neither
of these can occur, because $\sqrt{2}$ is irrational.
\end{proof}

If $p$ is a prime for which the real number $\sqrt{p}$ is normal
in base~$p$ --- a condition one expects every prime to satisfy ---
then a similar argument shows that there are
infinitely many exceptional powers of~$p$.  

\section{Improved bounds}
\label{S-improvements}

In this section we will explain how we obtained the improvements listed 
in Tables~\ref{Tbl-char2} and~\ref{Tbl-char3}.
Some of the entries in the table are immediate consequence of the 
corollaries in Section~\ref{S-double}, but other entries require a 
more detailed analysis.  We have written a Magma program that
carries out some of this analysis for us; 
let us begin by explaining what the program does.

Given a prime power $q$ and two positive integers $g$ and $N$,
we want to determine, if we can, whether there exists a
curve of genus $g$ over $\F_q$ with $N$ rational points.
Our program enumerates all of the polynomials that might 
possibly be the real Weil polynomial $h$ of such a curve,
where by ``might possibly'' we mean that
\begin{itemize}
\item all of the roots of $h$ are real numbers in the 
      interval $[-2\sqrt{q}, 2\sqrt{q}]$, and
\item the number of places of degree $d$ (for $d=1,\ldots,g$)
      predicted by $h$ are non-negative and in accord with the 
      Weil-Serre bounds.
\end{itemize}
The enumeration is carried out in one of two ways:  If the value of
$N$ corresponds to a defect of $6$ or less, the program uses precomputed 
tables of totally positive polynomials of deficiency at most $6$ 
(calculated as in~\cite{Smyth}) to list all of the appropriate
polynomials.  Otherwise, the program uses the algorithm 
from~\cite{Lauter:conf} to compute the appropriate~$h$'s.
For each candidate $h$, the program then uses the criterion 
of~\cite{Tate} to determine whether $h$ actually is the real Weil 
polynomial of an isogeny class of abelian varieties over $\F_q$.
For each $h$ that is a real Weil polynomial,
the program uses the factorization of $h$ to loop  through all 
of the splittings of the associated isogeny class into a product 
of lower-dimensional isogeny classes.  
It then computes the value of $r(F+V)$ for each such splitting.
If $r(F+V)=1$ then we know from Theorem~\ref{T-resultant}(a) that there is no 
curve with $h$ as its real Weil polynomial.  If $r(F+V)=2$ then the program tries 
to use Theorem~\ref{T-resultant}(b) and Lemma~\ref{L-obvious} to show that 
there is no curve with the given real Weil polynomial. 
(The program only uses Lemma~\ref{L-obvious}(a) with $d=1$,
and it only uses Lemma~\ref{L-obvious}(c) with all of the $d_i$
less than or equal to $g$.) 
If one of the isogeny classes in the splitting is an isogeny
class of elliptic curves and if the conclusion of Proposition~\ref{P-ec} 
is not satisfied, then again we know that there is no curve with the
given real Weil polynomial. If a polynomial $h$ is not 
eliminated by these filters, the program flags it as such and prints out 
three items:
\begin{itemize}
\item[(1)] The number of places of degree $d$ (for $d=1,\ldots,g$)
           that a curve would have to have in order to have $h$ 
           for its real Weil polynomial,
\item[(2)] the factorization of $h$, and
\item[(3)] a matrix giving the resultants of each pair of prime
           factors of $h$.
\end{itemize}
Likewise, if a polynomial $h$ {\em is\/} eliminated, the program will
print out item (2) above, together with an explanation of why
it eliminated the polynomial.  The program will also print out 
item (1) if it had to calculate that information in order to 
eliminate the polynomial.

For some specific choices of $q$, $g$, and $N$, our program eliminates
all possible real Weil polynomials.  For other choices there are only a few
real Weil polynomials left to consider, and sometimes we can eliminate these
by other methods; see Sections~\ref{S-descent}, \ref{S-exhaustion}, 
and~\ref{S-Hermitian} for examples of some of these methods.

Throughout this section, the symbol $m$ will always stand for the
integer $[2\sqrt{q}]$, where $q$ is the prime power currently
under discussion.

It was proven in~\cite{Serre:notes} 
(see also~\cite[Prop.~2]{LauterSerre:JAG})
that defect-$1$ curves are never possible when the genus is 
bigger than~$2$.  We will frequently use this fact without comment.

\subsection{Improvements for $q=4$}
\label{SS-4}
\subsubsection*{The case $q=4$, $g=5$, $N=18$}
We ran our Magma program for the case $q=4$, $g=5$, $N=18$.
The output is reproduced in the first Appendix.
The program finds eight polynomials
$h$ that might possibly be real Weil polynomials for a genus-$5$ curve over
$\F_4$ with $18$ points.
The first of the eight possibilities turns out not to be a real Weil polynomial;
it fails the local criterion given in~\cite{Tate}.

The second, fourth, fifth, seventh, and eighth possibilities are eliminated
by Theorem~\ref{T-resultant}(a).  For example, for the fifth possibility,
the program finds that $h$ factors
as $(x+2)^2 (x+4) (x^2 + 5x + 5)$.  If we let $h_1$ be the product of
the first and second of these factors (as the line
\begin{verbatim}
Splitting = [ 1, 2 ]
\end{verbatim}
in the output indicates we should do) and if we let $h_2$ be the
third factor, then the resultant of the radical of $h_1$ and the radical
of $h_2$ is $1$.

The third and sixth possibilities are eliminated by Theorem~\ref{T-resultant}(b)
and Lemma~\ref{L-obvious}.  For example, the polynomial $h$ for the
third possibility is $(x+1)(x+2)(x+3)^2(x+4)$.  We can factor this 
as $h_1 h_2$ where $h_1 = (x+2)$ and $h_2 = (x+1)(x+3)^2(x+4)$,
and then the resultant of the radicals of $h_1$ and $h_2$ is $2$.
Thus Theorem~\ref{T-resultant}(b) shows that any curve 
with $h$ as its real Weil polynomial 
would have to be a double cover of a curve $D$ whose real Weil polynomial
is either $h_1$ or $h_2$.
But if $D$ had $h_1$ as its real Weil polynomial then it would have $7$ points,
and this would contradict Lemma~\ref{L-obvious}(a) with $d=1$; while if
$D$ had $h_2$ as its real Weil polynomial then it would have genus~$4$,
and this would contradict Lemma~\ref{L-obvious}(b).
This argument is summarized in the lines
\begin{verbatim}
Splitting = [ 2 ]
Reasons: point counts, Riemann-Hurwitz
\end{verbatim}
of the output.

So we see that there is no genus $5$ curve over $\F_4$ with $18$ points.  
A curve with $17$ points is known, so we obtain the first entry
in Table~\ref{Tbl-exact}.

\subsubsection*{The cases $q=4$, $g\in\{10, 11\}$}
The improvements listed in Table~\ref{Tbl-char2} for
$q=4$ and $g\in\{10,11\}$ come from running our program.
Note that a curve over $\F_4$ of genus $10$ with $27$ points
is known, so we get the second entry in Table~\ref{Tbl-exact}.

\subsection{Improvements for $q=8$}
\label{SS-8}

\subsubsection*{The cases $q=8$, $g=5$, $N=32$ and $31$}
Our program shows that no genus-$5$ curve over $\F_8$ can
have exactly $32$ points.

The case $N=31$ is more interesting.  Our program shows that if
$C$ is a genus-$5$ curve over $\F_8$ with $31$ points, then
its real Weil polynomial must be 
$$h = (x+5)^3 (x^2 + 7x + 8).$$
The resultant of $x+5$ and $x^2 + 7x + 8$ is $2$, so $C$ is a double
cover of a curve $D$ whose real Weil polynomial is either
$(x+5)^3$ or $x^2 + 7x + 8$.  If $D$ had $(x+5)^3$ for its real Weil 
polynomial then $D$ would have genus~$3$, and this would
contradict Lemma~\ref{L-obvious}(c) since $C$ has an odd 
number of rational points.
Thus the real Weil polynomial of $D$ must be $x^2 + 7x + 8$.

We see that the Weil polynomial of $C$ must be
$$f_C = (x^2 + 5x + 8)^3 (x^4 + 7x^3 + 24x^2 + 56x + 64)$$
and the Weil polynomial of $D$ must be
$$f_D = x^4 + 7x^3 + 24x^2 + 56x + 64.$$
In Section~\ref{S-descent} we will use a 
Galois descent argument to show that this cannot occur.

\subsubsection*{The cases $q=8$, $g=9$, $N=47$ and $46$}
Our program shows that a genus-$9$ curve over $\F_8$ cannot have 
exactly $47$ points, and that if such a curve has exactly $46$ points
then its real Weil polynomial is either
$(x+3)(x+4)^3(x+5)^3(x^2+7x+9)$ or $(x+3)^4 (x+5)^5$.  We will
consider each of these two possibilities in turn.

Suppose $C$ is a genus-$9$ curve over $\F_8$ with real Weil 
polynomial $$h=(x+3)(x+4)^3(x+5)^3(x^2+7x+9).$$   
Since the resultant of $x+5$ and
$(x+3)(x+4)(x^2+7x+9)$ is $2$, the curve $C$ must be a double cover
of a curve $D$ over $\F_8$ whose real Weil polynomial is either
$(x+5)^3$ or $(x+3)(x+4)^3(x^2+7x+9)$.  The second option would 
require $D$ to have larger genus than is allowed by
Lemma~\ref{L-obvious}(b), so $D$ must have
real Weil polynomial $(x+5)^3$.  In particular, $D$ must have
exactly $24$ places of degree~$1$.

We see from the real Weil polynomial of $C$ that $C$ has no places 
of degree~$2$.  In particular, no rational places of $D$ can be inert 
in the double cover $C\to D$.  Since $C$ has $46$ rational places,
it must be the case that $22$ of the $24$ rational places of $D$ split
in the cover $C\to D$ and the other $2$ rational places ramify.

We also see from the real Weil polynomial of $C$ that $C$ has 
$109$ places of degree~$3$.
As we argued in the proof of Lemma~\ref{L-obvious}, the fact that $C$
has an odd number of degree-$3$ places implies that at least one
degree-$3$ place of $D$ ramifies.  Thus, at least $5$ geometric points
of $D$ ramify in the double cover $C\to D$.  Since the ramification
is necessarily wild, each ramification point contributes at least $2$
to the degree of the different of the cover, which means that the
degree of the different is at least~$10$.  But the Riemann-Hurwitz formula
shows that the degree of the different of the double cover
$C\to D$ is equal to~$8$.   
This contradiction shows that there is no genus-$9$ curve over $\F_8$ with
$(x+3)(x+4)^3(x+5)^3(x^2+7x+9)$ for its real Weil polynomial.

Suppose $C$ is a genus-$9$ curve over $\F_8$ with real Weil 
polynomial $(x+3)^4 (x+5)^5.$   Since the resultant of $x+3$ and
$x+5$ is $2$, the curve $C$ must be a double cover of a curve $D$
whose real Weil polynomial is either $(x+3)^4$ or $(x+5)^5$.
But the first option is impossible, because in that case
$D$ would have only $21$ rational points, which
contradicts Lemma~\ref{L-obvious}(a).
The second option is impossible as well, because in that case the 
genus-$5$ curve $D$ would have $34$ rational points, 
whereas we know that $N_8(5) \le 30$.

Thus there are no genus-$9$ curves over $\F_8$ with $46$ points.
A curve with $45$ points is known, so we get the third
entry in Table~\ref{Tbl-exact}.

\subsubsection*{The cases $q=8$, $g\in \{7,8,10,11,15\}$.}
The improvements we get when $q=8$ and $g\in \{7,8,10,11,15\}$
can all be obtained by running our program.

\subsection{Improvements for $q=16$}
\label{SS-16}

\subsubsection*{The cases $q=16$, $g\in \{4,5,7\}$.}
These improvements come directly from Theorem~\ref{T-square}(b).
Since a genus-$4$ curve over $\F_{16}$ with $45$ points is known,
we obtain the fourth entry in Table~\ref{Tbl-exact}.

\subsubsection*{The cases $q=16$, $g\in \{8,11,13,14\}$.}
The improvements we list in these cases are all obtained by running
our program.

\subsection{Improvements for $q=32$}
\label{SS-32}

\subsubsection*{The case $q=32$, $g=4$, $N=75$}
Suppose $C$ is a genus-$4$ curve over $\F_{32}$ with exactly $75$ 
rational points.  Then $C$ has defect $2$, so it must be of
one of the seven types listed in~\cite{LauterSerre:JAG}.
The type $[m,m,m-1,m-1]$ is forbidden by Theorem~\ref{T-resultant}(a),
and the fractional part of $2\sqrt{32}$ is small enough to
eliminate five of the others.  Thus $C$ must have 
type $[m,m,m,m-2]$, where $m=[2\sqrt{32}]=11$.
Theorem~\ref{T-resultant} tells us that $C$ 
must be a double cover of either a genus-$3$ curve (which is impossible, 
by Lemma~\ref{L-obvious}(b)) or of an elliptic curve whose Weil polynomial 
is $x^2 + (m-2)x + 32 = x^2 + 9x + 32$.  In Section~\ref{SS-32.4.75} we will 
show how the set of all genus-$4$ double covers of such elliptic curves 
can be enumerated.  We will see that none of the curves has $75$ points.

\subsubsection*{The cases $q=32$, $5\le g \le 15$.}
For $q=32$ and $g \ge 3$, a Galois descent 
argument~\cite{LauterSerre:JAG} shows that the 
Weil-Serre upper bound cannot be met,
and the previously-known best upper bound for $3 \le g \le 15$ 
was $q+1+gm-2$.  As we saw above, the only possible defect-$2$
zeta function is of type $[m,\ldots,m,m-2]$.
However, Corollary~\ref{C-type2} rules out this type of zeta function
when $g\ge 5$, so defect $2$ is impossible when $g\ge 5$.

Likewise, the arguments from the appendix of~\cite{Savitt} show
that for $g\ge 9$ the only possible defect-$3$ zeta function
for $q = 32$ is of type $[m,\ldots,m,m-3]$.  (This also depends on
the fact that the fractional part of $2\sqrt{32}$
is relatively small.)  
But Corollary~\ref{C-ec} shows that then $g\le 8$,
so defect $3$ is impossible when $g\ge 9$.
Thus our new upper bound is $q+1+gm-3$ for $5 \le g\le 8$,
and is $q+1+gm-4$ for $9\le g\le 15$.

\subsection{Improvements for $q=64$}
\label{SS-64}
\subsubsection*{The cases $q=64$, $11\le g \le 27$, $g\neq 12$}
If $g=11$, then  a curve meeting the Weil-Serre bound is not possible 
due to the results of Korchmaros-Torres~\cite{KorchmarosTorres}. 
Defect $2$ is also impossible, by Corollary~\ref{C-defect2}. 

We know from~\cite{FuhrmannTorres} that 
there is also no defect-$0$ curve when $13 \le g \le 27$, and 
it was shown in~\cite{LauterSerre:JAG} that defect $2$ is ruled out
by the Honda-Tate theorem.
But Theorem~\ref{T-square} shows that defects $3$ and 
$4$ are not possible either, so we get an upper bound of~$q+1+gm-5$.

\subsection{Improvements for $q=128$}
\label{SS-128}

Apart from the case $g=9$ and $N=324$ (explained below), all of our
improved bounds for $q=128$ can be obtained by running
our program.  However, we will take some time here to indicate how the
structure apparent in the $q=128$ results is a consequence of the
fact that $128$ is exceptional (in the terminology of
Section~\ref{S-exceptional}).
To simplify our discussion, let us introduce some terminology.

Suppose $h$ is a monic irreducible polynomial in $\Z[x]$, all
of whose roots in $\C$ are real and have magnitude at 
most $2\sqrt{q}$.  By the Honda-Tate theorem 
there is an integer $e>0$ such
that a power $h^n$ of $h$ is the real Weil polynomial of an abelian variety
over $\F_q$ if and only if $n$ is divisible by $e$.
We will say that $h^e$ is an {\em elementary\/} 
real Weil polynomial.  For example, the polynomial $(x+22)^7$
is an elementary real Weil polynomial over $\F_{128}$.

We define the {\em defect\/} of a real Weil
polynomial $h$ over $\F_q$ to be $m\deg h + \trace h$,
where the trace of a polynomial is as defined in
Section~\ref{S-exceptional}.  Note that if $C$ is a curve
over $\F_q$ of defect $d$ then its real Weil polynomial has
defect~$d$.  Also, the defect of a product of real Weil polynomials is
the sum of the defects.

Suppose $h\in\Z[x]$ is the real Weil polynomial of a curve $C$ 
over $\F_q$.
Let $H(x) = h(x - m - 1)$, so that all of the roots
of $H$ are positive real numbers.  
One checks that the defect of $h$ is the deficiency of $H$,
as defined in Section~\ref{S-exceptional}.
Smyth~\cite{Smyth} has written down all irreducible
monic polynomials $H$ in $\Z[x]$ with totally positive
roots and with deficiency at most $6$,
and using Smyth's work and the Honda-Tate theorem it is
not hard to write down a list of all of
the elementary real Weil polynomials $h$ over $\F_q$
with defect at most~$6$.  (A Magma program to reproduce Smyth's
work is available at the URL mentioned in the acknowledgments.)

There is only one elementary real Weil polynomial of defect $0$,
namely $(x + m)^\delta$, where $\delta$ is the defect-$0$ dimension of $q$.
Let us say that a real Weil polynomial over $\F_q$ is
{\em minimal\/} if it is coprime to $x + m$.
Given the list of elementary real Weil polynomials over $F_q$ of
defect at most~$6$, it is a simple matter to make a list
of all of the minimal real Weil polynomials over $\F_q$ of defect at most $6$.

Now suppose one is interested in genus-$g$ curves $C$ over $\F_q$ with
defect~$d\le 6$.  The real Weil polynomial of $C$ must be of the
form $(x + m)^n h$, where $h$ is a minimal real Weil polynomial
of defect $d$.  As we just noted, one can easily list these polynomials;
the task is made even simpler by the fact that only $h$ of certain
degrees can occur, since $n = g - \deg h$ must be a multiple
of the defect-$0$ dimension of $q$.  Furthermore, one
can use Theorem~\ref{T-resultant}(a) to 
exclude certain polynomials $(x + m)^n h$.

For instance, consider the case where $q=128$ and $g \equiv 2 \bmod 7$,
with~$g>2$.
There can be no defect-$0$ curves of genus $g$ because 
the defect-$0$ dimension of $q$ is $7$. 
There are no defect-$1$ real Weil polynomials because we took~$g>2$.
The only possible defect-$2$ polynomials are
$y^{g-2} (y-1)^2$ and $y^{g-2} (y^2 - 2y - 1)$, where $y=x+m$, but these
are eliminated by Theorem~\ref{T-resultant}(a).  
The possible defect-$3$ polynomials are
$y^{g-2} (y^2 - 3y + 1)$ and $y^{g-2} (y^2 - 3y - 1)$ and
$y^{g-2} (y^2 - 3y - 2)$.  The first two are eliminated by
Theorem~\ref{T-resultant}(a), and when $g > 9$ the third is eliminated
by Theorem~\ref{T-resultant}(b) and Lemma~\ref{L-obvious}.
The possible defect-$4$ polynomials are
$y^{g-2}(y^2 - 4y - 1)$ and $y^{g-2}(y^2 - 4y + 1)$ and 
$y^{g-2} (y-1)(y-3)$.  The first two are eliminated by
Theorem~\ref{T-resultant}(a), and when $g \ge 9$ the third
is eliminated by Theorem~\ref{T-resultant}(b) and Lemma~\ref{L-obvious}.
For defect $5$ there are several possible polynomials that
we cannot eliminate using our theorems.  
Combining all of the above, we see that when $g > 9$ is
congruent to $2\bmod 7$, we have
$N_{128}(g) \le q + 1 + mg - 5$.
A similar analysis can be done for the other congruence classes modulo~$7$.

There is a known curve of genus $4$ over $\F_{128}$ with $215$
rational points, so we obtain the fifth entry in Table~\ref{Tbl-exact}.

\subsubsection*{The case $q=128$, $g=9$, $323\le N \le 327$.}
The defect-$0$ dimension of $q$ is $7$, so there is no defect-$0$
curve of genus~$9$.  Defect $1$ is impossible because $g>2$.
The cases $N=325$ and $N=323$ are eliminated by our program.  The only
case remaining is $N=324$.

Our program shows that a genus-$9$ curve $C$ over $\F_{128}$ with $324$
points would have to have real Weil polynomial $(x+22)^7(x^2+41x+416)$
and would have to be a double cover of a genus-$2$ curve $D$  with real Weil
polynomial $x^2 + 41x + 416$.  From their real Weil polynomials, we see
that $C$ and $D$ each have $2$-rank $1$; that is, the $\F_2$-dimension of
the geometric $2$-torsion of their Jacobians is~$1$.  But then
the Deuring-Shafarevich formula (see~\cite{Subrao} and the references
listed in~\cite[\S 3]{Bouw}) shows that the double cover $C\to D$
must be unramified, which is clearly impossible.

\subsection{Improvements for $q=3$}
\label{SS-3}

\subsubsection*{The case $q=3$, $g=6$, $N=15$}
Running our program on this case leaves us with three real Weil polynomials
to consider.

The first is $(x+2)^2(x+3)(x^3+4x^2+x-3)$.  Factoring this as
$(x+3)$ times $(x+2)^2(x^3+4x^2+x-3)$ and applying Proposition~\ref{P-ec},
we find that a curve with this real Weil polynomial must be a triple
cover of an elliptic curve with Weil polynomial $x^2 + 3x + 3$.
We will show in the second Appendix that no such 
triple cover can have $15$ points.

The second real Weil polynomial we must consider is 
$(x+2)^2(x^2 + 3x - 1)(x^2 + 4x + 2)$.  Factoring this as
$(x^2 + 4x + 2)$ times $(x+2)^2(x^2 + 3x - 1)$ and applying 
Theorem~\ref{T-resultant}(b), we find that a curve with this
real Weil polynomial must be a double cover of a genus-$2$ curve
with real Weil polynomial $(x^2 + 4x + 2)$.  Searching through 
the genus-$2$ curves over $\F_3$, we find that there is exactly
one curve with that real Weil polynomial; it is given by the 
equation $y^2 = x^6 + x^5 + x^4 + x^2 - x + 1$.  In 
Section~\ref{SS-3.6.15.double} we will show that there is no
genus-$6$ double cover of this curve with $15$ points.

The third real Weil polynomial we are left to consider is
$(x+1)^2(x+3)^2(x^2+3x-1)$.  Writing this polynomial as the product
of $(x+3)^2$ and $(x+1)^2(x^2+3x-1)$ and applying Theorem~\ref{T-resultant}(b),
we find that a curve $C$ with this real Weil polynomial must be a 
double cover of a genus-$2$ curve with real Weil polynomial
$(x+3)^2$.  Such a genus-$2$ curve would have $10$ rational points --- but
this is impossible, because $N_3(2) = 8$.

Thus there is no genus-$6$ curve over $\F_3$ with $15$ rational points.
A curve with $14$ points is known, so we obtain the sixth entry 
in Table~\ref{Tbl-exact}.

\subsection{Improvements for $q=9$}
\label{SS-9}

\subsubsection*{The case $q=9$, $g = 13$, $N=66$}
Running our program shows that the only possible real Weil
polynomial in this case is $(x + 2)(x + 4)^6(x + 5)^6$.
Writing this polynomial as the product of $(x+4)^6$
with $(x+2)(x+5)^6$ and applying Theorem~\ref{T-resultant}(b),
we see that a genus-$13$ curve $C$ over $\F_9$ with $66$ points
must be a double cover of a curve $D$ such that either 
\begin{enumerate} 
\item the curve $D$ has $34$ rational points and has genus $6$,
and the double cover $C\to D$ is ramified at $4$ geometric points, or
\item the curve $D$ has $42$ rational points and has genus $7$, and
the double cover $C\to D$ is unramified.
\end{enumerate}
We note that the real Weil polynomial of $C$ shows that it has
no places of degree~$2$, so that no rational point of $D$ can
be inert in the double cover $C\to D$, and so that every 
degree-$2$ place of $D$ must be inert in $C\to D$.

Suppose $D$ has genus $6$.  Since $D$ has $34$ rational points
and none of them are inert in $C\to D$, and since $C$ has $66$ rational
points, we see that exactly $2$ rational points of $D$ are ramified.
Since there are $4$ geometric ramification points, a degree-$2$ place
of $D$ must ramify as well --- but we have just seen that every 
degree-$2$ place of $D$ must be inert, a contradiction.

Suppose $D$ has genus $7$.  Since no rational point of $D$ can be
inert or ramified in the double cover $C\to D$, each of the $42$
rational points of $D$ must split.  But then $C$ would have to have
$84$ rational points, contradicting the fact that it has only~$66$.  

Thus neither of the two possibilities listed above can hold, and
there can be no genus-$13$ curve over $\F_9$ with $66$ points.

\subsubsection*{The cases $q=9$, $g \in\{9,10,11,12,14,15,16,17,18\}$}
The improvements listed in Table~\ref{Tbl-char3} for $q = 9$
and $g \in\{9,10,11,12,14,15,16,17,18\}$ can all be obtained
simply by running our program.  Note that there is a known
genus-$10$ curve over $\F_9$ with $54$ points,
so we get the seventh entry in Table~\ref{Tbl-exact}.

\subsection{Improvements for $q=27$}
\label{SS-27}

\subsubsection*{The cases $q=27$, $g=4$, $N=66$ and $N=65$}
Suppose $C$ is a genus-$4$ curve over $\F_{27}$ with exactly $66$ 
rational points.  Then $C$ has defect $2$, and of the seven
types of zeta function from~\cite{LauterSerre:JAG} 
one is eliminated by Theorem~\ref{T-resultant}(a) and 
five more are forbidden by the size of the fractional part
of $2\sqrt{27}$.  The only possibility remaining is $[m,m,m,m-2]$,
where $m=[2\sqrt{27}]=10$. Theorem~\ref{T-resultant}(b) tells us that $C$ 
must be a double cover of either a genus-$3$ curve (which is impossible, 
by Lemma~\ref{L-obvious}(b)) or of an elliptic curve whose Weil polynomial 
is $x^2 + (m-2)x + 27 = x^2 + 8x + 27$.  In Section~\ref{SS-27.4.66} we will 
show how the set of all genus-$4$ double covers of such elliptic curves 
can be enumerated.  We will see that none of the curves has $66$ points.

We show in the second Appendix that there is no 
genus-$4$ curve over $\F_{27}$ with exactly $65$ 
rational points.

A genus-$4$ curve over $\F_{27}$ with $64$ points is known, so we 
obtain the eighth entry in Table~\ref{Tbl-exact}.

\subsubsection*{The cases $q=27$, $5 \le g \le 13$}
First we note that~\cite[Thm.~1]{LauterSerre:JAG}
shows that the Weil-Serre bound cannot be met when~$g \ge 3$.  

Now we show that defect $2$ is impossible for $g\ge 5$.
For $g>5$ this follows from Corollary~\ref{C-defect2}.
When $g=5$ Corollary~\ref{C-type2} shows that $[m,m,m,m,m-2]$
is not a possible type, and the proof of Corollary~\ref{C-defect2}
shows that the only other possible type is 
$[m,m,m,m+\sqrt{3}-1,m-\sqrt{3}-1]$.  But this last type is
also impossible, because the fractional part of $2\sqrt{27}$ is 
less than $\sqrt{3}-1$.

Finally we note that the appendix to~\cite{Savitt} shows that
the only defect-$3$ curves over $\F_{27}$ when $g\ge 9$ are
of type $[m,\ldots,m,m-3]$, and Corollary~\ref{C-ec} shows that 
this type is impossible for $g\ge 9$.

\subsubsection*{The case $q=27$, $g=14$, $N=164$}
Our program eliminates this possibility.

\subsection{Improvements for $q=81$}
\label{SS-81}

\subsubsection*{The cases $q=81$, $13 \le g \le 35$, $g \ne 16$}
{}From \cite{FuhrmannTorres} and \cite{KorchmarosTorres}, 
we know that no defect-$0$ curves are possible 
for $13 \le g \le 35$, $g \ne 16$.  But defect $2$ and $3$ are not 
possible either by Theorem~\ref{T-square}, so the
upper bound for these cases is at most~$q+1+gm-4$.
When $g\ge 18$ we see from Theorem~\ref{T-square}(c)
that defect $4$ is impossible as well, so our new upper bound
for $18\le g \le 35$ is~$q+1+gm-5$.

\subsection{Cases where few Weil polynomials are possible}
\label{SS-no}

We have tried to use our program to obtain further improvements to
the upper bounds listed in the van der Geer-van der Vlugt tables,
but it appears that we have already picked most of the low-hanging fruit.
For example, for every $q$ listed in the tables, and for every $g\le 10$,
we have taken the best current upper bound $N$ for $N_q(g)$ and run
our program on the triple $(q,g,N)$.    In each such case, our program
indicates that there are real Weil polynomials that are not eliminated
by the criteria that we built in to the program.  This is not to say,
however, that our methods cannot give further improvements in these cases:
For instance, for $q = 27$ and $5\le g\le 8$, a curve meeting the
best current upper bound would have to be a triple cover of a defect-$3$
elliptic curve, and it might be possible to use the description of 
such covers given in Section~\ref{S-char3} to enumerate all of the
possible curves.

We mention just two more interesting cases.  The smallest genus $g$
for which $N_2(g)$ is not known is $g = 12$; it is known that 
$N_2(12)$ is either $14$ or $15$.  Running our program on the
case $q = 2$, $g=12$, $N=15$ took almost $18$ hours using Magma~$2.8$
on a $2$ GHz Pentium~$4$, 
and we found that there are eight possible
real Weil polynomials to consider.  We were unable to eliminate all of
these polynomials, and we were not able to use them to direct a search for
a genus-$12$ curve with $15$ points.

On the other hand, running our program on the case $q = 4$, $g = 7$, $N=22$
produces six candidate real Weil polynomials, and by a number of
{\em ad hoc\/} arguments we were able to eliminate all but one of them
from consideration.  We find that if a genus-$7$ curve over $\F_4$ has
$22$ points, then its real Weil polynomial must be $x(x+2)^2(x+3)^3(x+4)$.
Arguments along the lines of those provided in the second Appendix can be
used to eliminate this possibility as well; we will provide details in
a forthcoming paper.

\section{A Galois descent argument}
\label{S-descent}

In Section~\ref{SS-8} we showed that a genus-$5$ curve $C$ over $\F_8$
having exactly $31$ points must be a double cover of a genus-$2$ curve~$D$. 
In this section we will use a Galois descent argument to show that the 
curves $C$ and $D$ and the degree-$2$ map $C\to D$ can all be defined 
over $\F_2$, and we will show how this leads to a contradiction.

Let $f_2 = x^4 + x^3 + 2x + 4$ and let $g_2 = x^2 - x + 2$.
Let $\pi$ be a root of $f_2$ in $\Qb$ and $\rho$ be a root of $g_2$ in $\Qb$.
Let 
$$f_8 = x^4 + 7x^3 + 24x^2 + 56x + 64 
  \quad\text{and}\quad 
  g_8 = x^2 + 5x + 8.$$
Note that $\pi^3$ is a root of $f_8$ and that $\rho^3$ is a root of $g_8$. 
The arguments from Section~\ref{SS-8} show that it will suffice for us to 
prove the following:

\begin{prop}
\label{P-8-5-31}
There is no genus-$5$ curve over $\F_8$ with Weil 
polynomial~$f_8 g_8^3$.
\end{prop}

\begin{proof}
We know from Section~\ref{SS-8} that any such curve must be a double 
cover of a genus-$2$ curve $D$ with Weil polynomial~$f_8$. 
Let us first identify the curve~$D$.

\begin{claim}
There is exactly one principally polarized abelian surface 
over $\F_8$ with Weil polynomial equal to~$f_8$.  
It is the polarized Jacobian of the curve $y^2 + xy = x^5 + x$. 
\end{claim}

\begin{proof}
Every such principally-polarized variety is a Jacobian,
because the varieties in the isogeny class determined by $f_8$ 
are absolutely simple (see~\cite[Thm.~6]{HoweZhu}). By explicitly
enumerating the genus-$2$ curves over $\F_8$ one finds that the
curve given above is the only curve whose Jacobian has 
Weil polynomial~$f_8$.
\end{proof}

Suppose, to get a contradiction, that $C$ is a genus-$5$ curve
over $\F_8$ with Weil polynomial $f_8 g_8^3$.  Since the resultant
of the real Weil polynomials associated with $f_8$ and $g_8$ is $2$,
Lemma~\ref{L-resultant} shows that there is an exact sequence
$$0 \to \Delta \to A \times B \to \Jac C \to 0,$$ where 
$A$ and $B$ are abelian varieties over $\F_8$ with Weil polynomials
$f_8$ and $g_8^3$, respectively, and where the projections
$A \times B \to A$ and $A \times B \to B$ induce monomorphisms 
$\Delta \hookrightarrow A[2]$ and
$\Delta \hookrightarrow B[2]$.
Since $\Jac C$ is a Jacobian and hence has a principal polarization,
Lemma~\ref{L-resultant} shows that $\Delta$ is self-dual.
Furthermore, $\Delta$ is nontrivial, because the principal
polarization on $\Jac C$ is indecomposable.

Every finite group-scheme $G$ in characteristic $p$ can be 
written as a product of four sub-group-schemes:
$$G = \Grr\times\Grl\times\Glr\times\Gll,$$
where $\Grr$ is a reduced group-scheme whose Cartier dual 
is reduced, where $\Grl$ is a reduced group-scheme whose 
Cartier dual is local, and so on.  (See~\cite[\S I.2]{Oort}.)
A group-scheme of $p$-power rank in characteristic $p$ can 
have no reduced-reduced part. Furthermore, if $G$ is 
self-dual --- for example, if $G$ is the kernel of a 
polarization --- then $\Grl$ and $\Glr$ are\ dual to one another.

Now, $B$ is an ordinary abelian variety, and the kernel of 
multiplication-by-$p$ on an ordinary abelian variety in 
characteristic $p$ has no local-local part.  Thus $B[2]$ 
consists of a reduced-local factor of rank $8$ and a 
local-reduced factor of rank~$8$.

The variety $A$ is not ordinary, so $A[2]$ has a local-local 
component. However, $A$ has positive $2$-rank, so $A[2]$ has
a reduced-local component as well. The only possibility is 
that $A[2]$ has a reduced-local component of rank~$2$, a 
local-reduced component of rank~$2$, and a local-local 
component of rank~$4$.

There are supposed to be monic maps from $\Delta$ to $A[2]$ 
and to $B[2]$.  Since $\Delta$ can be viewed as a subscheme 
of $B[2]$ it can have no local-local part.  Thus, the 
monomorphism $\Delta\to A[2]$ must take $\Delta$ onto the 
product of the reduced-local and the local-reduced part of A[2]. 
Since $\Delta$ is self-dual, it follows that $\Delta$ has rank~$4$ 
and is the product of a rank-$2$ reduced-local group and a rank-$2$
local-reduced group.

As in the proof of Lemma~\ref{L-resultant}, let $\mu_A$ and 
$\mu_B$ be the degree-$4$ polarizations on $A$ and $B$ that
we get by pulling back the canonical polarization of $\Jac C$ 
via the map $A\times B\to \Jac C$.  We know that $\Delta$ is 
isomorphic to $\ker \mu_A$ and to $\ker \mu_B$. The local-reduced 
subgroup of $\ker \mu_A$ is maximal isotropic, so the polarization
$\mu_A$ on $A$ gives rise to a principal polarization on the 
quotient of $A$ by this subgroup. It follows from the claim 
we made above that this quotient variety is the Jacobian of 
the curve $D$.

We can make a diagram
\begin{equation*}
\begin{matrix}
A & \mapright{\mu_A} & \hat{A} \\
\Bdownarrow& & \Buparrow \\
\Jac D& \longrightarrow & \hat{\Jac D}
\end{matrix}
\end{equation*}
where the left arrow is the degree-$2$ isogeny $A \to \Jac D$, 
the right arrow is the dual of this isogeny, and 
the bottom arrow is the canonical polarization on~$\Jac D$.

Now, $\Jac D$ has exactly one reduced-local subgroup of order~$2$, 
and it is defined over $\F_2$.  It is in fact the kernel of 
multiplication by $1+\pi$, so $\Jac D$ divided by this subgroup 
is geometrically isomorphic to $\Jac D$. Now, the composition of the 
bottom and right arrows gives an isogeny $\Jac D \to\hat{A}$ whose 
kernel is reduced-local and of order 2. So geometrically, 
$\hat{A}$ is isomorphic to $\Jac D$, which means that $\hat{A}$ 
is a twist of~$\Jac D$.  But $\hat{A}$ is isogenous to $\Jac D$ 
over $\F_8$, and the quadratic twist of $\Jac D$ is not isogenous 
to $\Jac D$ (as we can see by checking Weil polynomials), 
so $\hat{A}$ must be isomorphic to $\Jac D$ over~$\F_8$.  
It follows that $A$ is isomorphic to $\Jac D$ as well.  Thus $A$, 
and the polarization $\mu_A$, can be defined over~$\F_2$.

On the other hand, $B$ and the polarization $\mu_B$ can be 
defined over $\F_2$ simply because 
$\Z[\rho,\bar{\rho}] = \Z[\rho^3,\rhobar^3]$. 
(This is essentially the Galois descent argument that Serre gives 
in~\cite{Serre:notes} and 
in the appendix to~\cite{LauterSerre:JAG}.) 
This means that the whole diagram 
\begin{equation*} 
\begin{matrix} 
A \times B  & \longrightarrow & \hat{A}\times\hat{B} \\ 
\Bdownarrow &                 & \Buparrow \\ 
\Jac C      & \longrightarrow & \hat{\Jac C} 
\end{matrix} 
\end{equation*} can be descended down to~$\F_2$. 

Now we want to know whether we can have a curve $C$ over $\F_2$ 
with Weil polynomial equal to~$f_2 g_2^3$. Again we find
that $C$ must be a double cover of a genus-2 curve~$D$ with 
Weil polynomial~$f_2$.  But then we find that $C$ has 
$13$ points over $\F_4$ and $D$ has $4$ points over~$\F_4$, and 
this is impossible. 
\end{proof}

\section{Exhaustive searches over small spaces} 
\label{S-exhaustion}

In this section we will give three examples that show how 
Theorem~\ref{T-resultant} can give us enough information about a 
curve with a certain number of points for us to have a computer 
look at every such curve and bound its number of rational points.

\subsection{The case $q=27$, $g=4$, $N=66$}
\label{SS-27.4.66}

We showed in Section~\ref{SS-27} that a genus-$4$ curve over $\F_{27}$ 
with exactly $66$ rational points must be a double cover of an elliptic 
curve with Weil polynomial $x^2 + 8x + 27$. 
There are exactly $4$ elliptic curves over $\F_{27}$ with this 
Weil polynomial;  one of them is defined over $\F_3$, and
the other three are Galois conjugates of one another.  Given such an 
elliptic curve $E$, we will show how the genus-$4$ double covers $C$ 
of $E$ can be enumerated by computer.

The function field of $C$ must be obtained from that of $E$ by adjoining
a root of $z^2 = f$, where $f$ is a function on $E$. By the 
Riemann-Hurwitz formula, in order for $C$ to have genus $4$ the 
divisor of $f$ must be of the form 
$$P_1 + \cdots + P_6 + 2D,$$ 
where the $P_i$ are distinct geometric points on $E$ and where $D$ 
is a divisor of degree~$-3$.  There is a function $g$ on $E$ such that
$$D + \divisor g = Q - 4\infty,$$ 
where $\infty$ is the infinite point on~$E$ and where $Q$ is a 
rational point on~$E$. Replacing $f$ with $fg^2$ does not change the
double cover of $E$. Thus, we may assume that $C$ is given by 
adjoining a root of $z^2 = f$, where $f$ is a function on $E$ 
whose divisor is of the form 
$$P_1 + \cdots + P_6 + 2Q - 8\infty.$$

We can also change the map $C\to E$ by following it with a translation
map on~$E$.  Translating  $E$ by a rational point $R$ has the effect 
of replacing $f$ with a function whose divisor is 
$$(P_1 + R) + \cdots + (P_6 + R) + 2(Q + R) - 8R$$ 
(where the sums in parentheses take place in the algebraic group~$E$). 
By modifying this new $f$ by the square of a function we can get the 
divisor of $f$ to be 
$$(P_1 + R) + \cdots + (P_6 + R) + 2(Q - 3R) - 8\infty.$$ 
If we choose representatives of the classes of $E(\F_{27})$ 
modulo $3E(\F_{27})$, then we may assume that $Q$ is one of these
representatives.  It turns out that for each of the possible curves 
$E$ the group $E(\F_{27})/ 3E(\F_{27})$ has order~$3$, so for each 
$E$ we need consider only $3$ possible $Q$'s.
We can choose our $Q$'s so that they do not lie in $E[2]$.

Let us write $E$ in standard Weierstrass form $y^2 = x^3 + ax^2 + bx + c$
and try to write down all of the functions $f$ as above in a standard form. 
There are two cases to consider, depending on whether or not
any of the $P_i$ is~$\infty$.

Suppose that one of the $P_i$ is~$\infty$.
Then $f$ has degree~$7$ and its only pole is at~$\infty$, 
and $f$ has a double zero at $Q$.  Since $Q$ is 
not a $2$-torsion point by assumption,
we may write $Q = (x_0,y_0)$ with $y_0\neq 0$. 
Note that then $x-x_0$ is a uniformizing parameter at~$Q$. 
Let $f_0$ be a linear polynomial that defines 
the tangent line to $E$ at~$Q$.
Then up to squares $f$ can be written as 
$$f = \pm(f_1 + c_0 f_0),$$
where $f_1$ is a function of the form
\begin{multline*}
(x-x_0)^2 \cdot \text{(polynomial in $x$ of degree $\le 1$)} \\
       + (x-x_0)(y-y_0) \cdot \text{(monic linear in $x$)}.
\end{multline*}
Likewise, if no $P_i$ is~$\infty$, then we may write     
$f = \pm(f_1 + c_0 f_0)$ where $f_1$ is a function of the form 
\begin{multline*}
(x-x_0)^2 \cdot \text{(monic quadratic in $x$)} \\
       + (x-x_0)(y-y_0) \cdot \text{(polynomial in $x$ of degree $\le 1$)}.
\end{multline*}
It is not hard at all to have a computer algebra system write 
down all of these possible $f$'s for a given $E$.

Now our problem is to count the points on the extension of 
$E$ defined by $z^2 = f$.  It is easy to get an overestimate: 
If $P$ is a rational point on $E$ for which $f(P)$ is a nonzero
square, then there are two rational points of $C$ lying above $P$. 
If $f(P)$ is not a square, then there are no rational points of 
$C$ above $P$.  If $P$ is a simple or a triple zero of $f$, 
then there is one rational point of $C$ above $P$.  And if 
$P$ is a double zero of $f$, then there are at most $2$ 
rational points of $C$ lying above~$P$.

What we actually did in practice for each candidate $f$ was to: 
\begin{enumerate} 
\item Eliminate $f$ from consideration if we could find more than three
      points $P$ on $E$ with $f(P)$ nonsquare;
\item Calculate the overestimate for $\#C(\F_{27})$ described above; 
\item Discard $f$ if the overestimate was less than $66$; 
\item Check to see that the divisor of $f$ was of the proper form,
      and discard $f$ if it was not.
\end{enumerate}

No candidate $f$'s made it through these filters, so we never had to 
worry about resolving the singularities of our model for $C$ to get 
an exact point count.

It took a little more than twelve hours using Magma $2.9$ 
on a $400$ MHz PowerPC G4 processor to search through all of the $(E,f)$ pairs that we 
had to consider. (Our Magma program is available at the URL 
mentioned in the acknowledgments.) Note that we need only consider 
two $E$'s; if one of the $E$'s that is defined only over $\F_{27}$ 
has a double cover with $66$ points, then so do all of its conjugates.

\subsection{The case $q=32$, $g=4$, $N=75$}
\label{SS-32.4.75}

We showed in Section~\ref{SS-32} that a genus-$4$ curve over $\F_{32}$ 
with exactly $75$ rational points must be a double cover of an elliptic 
curve with Weil polynomial $x^2 + 9x + 32$. 
There are exactly $5$ elliptic curves over $\F_{32}$ with this 
Weil polynomial, and they are all conjugate to one another 
over~$\F_2$.  (If $a\in \F_{32}$ satisfies $a^5 + a^2 + 1 = 0$ then 
the elliptic curve $E$ defined by 
$y^2 + xy = x^3 + x^2 + a^7$ has the correct Weil polynomial.) 
As in the preceding section, we can easily program a 
computer to enumerate the genus-$4$ double covers of such an 
elliptic curve and check to see whether any of these double covers
has $75$ points.  The only complication is that a double cover in 
characteristic $2$ is given by an Artin-Schreier extension of function
fields instead of a Kummer extension.  

Suppose $C$ is a double cover of the curve $E$ given above.  
Then the function field of $C$ is obtained from that of $E$ by adjoining 
a root of $z^2 + z = f$, where $f$ is a function on $E$.  The points of $E$
that ramify in the cover $C\to E$ are contained in the set of poles of~$f$;
to determine whether a pole $P$ of $f$ is a ramification point, and to 
determine the contribution of $P$ to the different of the extension $C\to E$,
we look at the expansion of $f$ in the local ring of $E$ at $P$. 
According to~\cite[Prop.~III.7.10]{Stichtenoth}, if there is a function 
$g_P$ such that $f + g_P^2 + g_P$ has no pole at $P$, then $P$ is unramified. 
If there is no such function, then we can at least find a function $g_P$ 
so that $f + g_P^2 + g_P$ has a pole of odd order at $P$.  If the pole has 
order $m$, then the differential exponent of $P$ in the extension $C\to E$ 
is $m+1$.

Suppose for each pole $P$ of $f$ we find a function $g_P$ as above. 
Then by Riemann-Roch we can find a function $g$ on $E$ that has 
poles only at $\infty$ and at the poles $P$ of $f$ and such that 
$g - g_P$ has no pole at $P$ for every $P\neq \infty$.  Replacing 
$f$ by $f + g^2 + g$ does not change the extension $C\to E$, but it 
allows us to assume that $f$ has only odd-order poles, except perhaps 
at infinity.  By modifying $f$ in this same way by functions with poles 
only at $\infty$, we may also assume that if $f$ has an even order pole 
at infinity, then the order of the pole is at most $2$.

Now suppose that $C$ has genus $4$ and has $75$ rational points. 
Then the Riemann-Hurwitz formula 
shows that there are three possible configurations for the different 
of $C\to E$:  There are either 
\begin{itemize} 
\item[(1)] three points with differential exponent $2$, 
\item[(2)] one point with differential exponent $2$ and one with differential
           exponent~$4$, or 
\item[(3)] one point with differential exponent $6$. 
\end{itemize} 
The second possibility cannot occur, because each of the ramification points
would have to be rational over $\F_{32}$, and this would force $C$ to 
have an even number of rational points.

Suppose we are in case (3).  Then the one ramification point $P$ is 
rational, and by following the map $C\to E$ with a translation by $-P$, 
we may assume that the point $P$ is the infinite point $\infty$ on $E$. 
Modifying the corresponding $f$ as above, we find that we may assume 
that $f$ is a function of degree~$5$ whose only pole is at~$\infty$. 
Thus we may assume that $f$ has the shape 
$$f = (ax + b) y + (cx + d)$$ 
where $a\neq 0$.  Furthermore, by modifying $f$ by constants of the 
form $e^2 + e$, we may assume that $d$ is either $0$ or $1$.

Suppose we are in case (1), with ramification at $P_1$, $P_2$, and $P_3$. 
Since $C$ has an odd number of rational points, at least one of the $P_i$ 
is rational.  If we label this point $P_3$ and then translate by $-P_3$, 
we find that we may assume that $P_3 = \infty$.  We may also assume 
that neither $P_1$ nor $P_2$ is the unique $2$-torsion point on $E$ 
(that is, the unique point with $x = 0$), because 
if (say) $P_1$ is the $2$-torsion point on $E$, we can translate by $-P_2$ 
so that the ramification locus becomes $\{P_1 - P_2, \infty, -P_2\}$. 
Thus we may write $P_1 = (x_1,y_1)$ and $P_2 = (x_2,y_2)$ with $x_1\neq 0$
and $x_2\neq 0$.  Then we may write $f$ in the form 
$$f = ax + b\frac{y + y_1 + x_1}{x + x_1} 
         + c\frac{y + y_2 + x_2}{x + x_2} + d$$
where $b$ and $c$ are nonzero, where $d$ is either $0$ or $1$, and where 
$a$ is nonzero if and only if $f$ has a pole of order $2$ at~$\infty$. 
Note that if $P_1$ and $P_2$ are not defined over $\F_{32}$ then they 
are quadratic conjugates of one another, and so $b$ and $c$ must be 
quadratic conjugates of one another in order for $f$ to be defined 
over~$\F_{32}$.

It is a simple matter to count points on the curve $C$ defined 
by $z^2 + z = f$, where $f$ is as above, because we are assuming 
that every pole of $f$ ramifies.  So 
if $P$ is a rational point on $E$ that is a pole of $f$,
then there is one point on $C$ lying above~$P$.
If $P$ is a rational point that is not a pole of~$f$, 
then there are either two or zero
rational points on $C$ over $P$, depending on whether the
trace of $f(P)$ to $\F_2$ is $0$ or~$1$.

We used Magma to enumerate all of the 
possible $f$'s for one of the elliptic curves $E$ given above. 
(Our Magma program is available at the URL 
mentioned in the acknowledgments.) 
For each $f$ we counted points on the curve $z^2 + z = f$.
No $f$ gave us $75$ points.
Thus we verified that there is no genus-$4$ curve over $\F_{32}$ 
having exactly $75$ rational points.

\subsection{The case $q=3$, $g=6$, $N=15$}
\label{SS-3.6.15.double}
In Section~\ref{SS-3} we showed that there were two possible real 
Weil polynomials for a genus-$6$ curve over $\F_3$ having $15$ points.
One of the two polynomials was $(x+2)^2(x^2+3x-1)(x^2+4x+2)$, and
we showed that any curve $C$ with this real Weil polynomial must be a 
double cover of the genus-$2$ curve $D$ defined by 
$y^2 = x^6 + x^5 + x^4 + x^2 - x + 1$.   We note that $C$ has $15$
places of degree~$1$ and $53$ places of degree~$5$, and the proof
of Lemma~\ref{L-obvious} shows that therefore a degree-$1$ place and
a degree-$5$ place of $D$ must ramify in the double cover $C\to D$.
(Since the degree of the different of the cover is $6$ by Riemann-Hurwitz,
no other places of $D$ can be ramified.)
In this section we will show how one can make a short list of 
double covers of $D$ that contains all of the genus-$6$ covers ramified
only at a degree-$1$ place and a degree-$5$ and having $15$ rational points.
We will find that there are no such double covers.

Note that the automorphism group of $D$ is cyclic of order~$8$, generated
by the map $(x,y) \mapsto \big((1+x)/(1-x), y/(1-x^3)\big)$.  This group
acts transitively on the rational points of $D$, so we may assume that
the rational ramification point of the double cover $C\to D$ is our
favorite rational point on~$D$.  We will choose this point to be
the rational point on $D$ that is a pole of
the function $x$ and a zero of the function $y - x^3$,
which point we will denote by $\infty^+$.

Let $K$ be the function field of $D$ and let $L$ be the function field of $C$.
Then there is a function $f$ on $D$ such that $L = K(z)$ for an element $z$
with $z^2 = f$.  Let $P$ be the degree-$5$ place at which $C\to D$ ramifies.
Then the divisor of $f$ is $P - 5\infty^+ + 2E$
for some degree-$0$ divisor $E$ on $D$.  By Riemann-Roch, there is 
a function $g$ on $D$ whose divisor is $F - 2\infty^+ - E$ for some
effective degree-$2$ divisor $F$.  Replacing $z$ with $zg$ and $f$ with $fg^2$,
we find that we may assume that the divisor of $f$ is $P + 2F - 9\infty^+$
for some effective divisor $F$ of degree~$2$.

The divisor $F$ must have one of four possible shapes, each
considered below.  For each possibility, we had Magma check that that
there is no function $f$ on $E$ whose divisor is of the right form
and that gives an extension with~$15$ points.  
(Our Magma routines for doing this
are available at the URL mentioned in the acknowledgments.)

\subsubsection*{Case {\rm 1:} $F$ consists of one place of degree~$2$}
Since $D$ has $8$ rational points and $C$ has $15$, we see that the rational
points of $D$ that do not ramify in the double cover $C\to D$ must split.
Since $f$ is nonzero at the rational points of $D$ other than $\infty^+$
and since these points all split, $f$ must evaluate to $1$ at these points.
It is a simple matter to enumerate all of the elements of the Riemann-Roch
space $\calL(9\infty^+)$ that evaluate to $1$ at the other rational
points of~$D$, and to check that none of them has a divisor of the form
$P + 2F - 9\infty^+$ for a degree-$5$ place~$P$.

\subsubsection*{Case {\rm 2:} $F$ consists of two possibly equal 
places of degree~$1$, neither equal to $\infty^+$}
To handle this case, we consider all possible pairs of points $F_1$ and $F_2$
on $D$.  For each pair, we have Magma enumerate the elements of 
$\calL(9\infty^+)$ that vanish at $F_1$ and $F_2$ and that evaluate to $1$
at the other rational points on $D$.  For each such function, we
check that its divisor is not of the form $P+2F_1 + 2F_2 - 9\infty^+$
for a degree-$5$ place~$P$.

\subsubsection*{Case {\rm 3:} $F$ consists of $\infty^+$ and some other
degree-$1$ place}
Now we loop over all rational points $F_1\neq\infty^+$ of $D$, and consider
the elements of $\calL(7\infty^+)$ that vanish at $F_1$ and that evaluate
to $1$ at the other rational points of~$D$.  For each such function, we
check that its divisor is not of the form $P+2F_1 - 7\infty^+$
for a degree-$5$ place~$P$.

\subsubsection*{Case {\rm 4:} $F$ consists of two copies of $\infty^+$}
For this case we must consider the elements of $\calL(5\infty^+)$ that evaluate
to $1$ at the other rational points of~$D$.  It turns out that the
only such function is the constant function~$1$.

Thus we find that there are no curves over $\F_3$ having real
Weil polynomial $(x+2)^2(x^2+3x-1)(x^2+4x+2)$.

\section{Triple covers of elliptic curves in characteristic $3$}
\label{S-char3}

Note: In the original version of this paper, the results of 
Section~\ref{SS-different} were incorrect; since we used these
results in Sections~7.3 and~7.4,
we have had to come up with new arguments for the results of those
sections.  These new arguments can be found in the second Appendix.

We therefore no longer have any need for the results of 
Sections~\ref{SS-standardform} and~\ref{SS-different}. 
Nevertheless, we have decided to include the corrected 
versions of these sections in this version of our paper, 
in case the results will be useful for some other purpose.
Sections~7.3 and~7.4 have been deleted.

\subsection{A convenient standard form}
\label{SS-standardform}
Suppose $k$ is a finite field of characteristic~$3$, suppose $E$ is an 
elliptic curve over $k$, and suppose $C$ is a curve over $k$ for which
there is a degree-$3$ map $C\to E$.  We will show that $C$ can be
given in a convenient standard form.  We will limit ourselves to 
covers $C\to E$ for which a certain assumption (stated below) holds.

Let $L$ and $K$ be the function fields for $C$ and $E$, respectively,
and view $L$ as a degree-$3$ extension of $K$ via the degree-$3$ 
map $C\to E$.
Choose a generator for $L$ over $K$ whose trace to $K$ is~$0$.
Then there are functions $f$ and $g$ in $K$ such that $z^3 - fz - g = 0$.
Suppose we write the divisor of $f$ in the form
$$\divisor f = P_1 + \cdots + P_n + 2D,$$
where the $P_i$ are distinct geometric points of $E$ (and where $n$ is
necessarily even).  Note that since $f$ is (up to squares) the
discriminant of the extension $L/K$, the number $n$ is the
number of points of $E$ at which the discriminant of $L/K$ has
odd valuation.

\begin{assumption}
We will assume that $n$ is coprime to $\#E(k)$.
\end{assumption}

Under this assumption, there is a rational point $Q$ on $E$
such that $nQ = P_1 + \cdots + P_n$ in the group of points of~$E$.  
By composing the given map $C\to E$ with a translation, 
we may assume that $Q$ is the infinite point $\infty$ on~$E$.  
By replacing the divisor $D$ above with $D + (n/2)\infty$, 
we may write
$$\divisor f = P_1 + \cdots + P_n - n\infty + 2D$$
where the $P_i$ are distinct geometric points on $E$ whose sum is $0$
in~$E$.  It follows that $D$ is a degree-$0$ divisor, and the sum (in $E$)
of the points in $D$ is a $k$-rational $2$-torsion point on~$E$.
But since $\#E(k)$ is coprime to the even number $n$, the
only $k$-rational $2$-torsion point on $E$ is~$\infty$.
Therefore $D$ is a principal divisor,
because it has degree~$0$ and the sum (in $E$) of its points is zero.
Write $D = \divisor h$ for some function $h$.  Replacing $z$, $f$, and $g$
with $z/h$, $f/h^2$, and $g/h^3$, respectively, we find that we still
have $z^3 - fz - g = 0$, but now the divisor of $f$ is
$$\divisor f = P_1 + \cdots + P_n - n\infty.$$
If one of the $P_i$ (say~$P_n$) is equal to~$\infty$,
replace $n$ by $n-1$ and delete the point~$P_n$ from the
expression for~$\divisor f$.
The integer $n$ may no longer be even, but now we 
have that the $P_i$ are all distinct and that none of them
is $\infty$.  We do at least know that $n$ is not~$1$, because 
there is no function on $E$ with divisor $P-\infty$.

Now suppose $P$ is a finite place of $E$ at which $g$ has a pole.
Suppose $\ord_P g$ is a multiple of $3$, say $\ord_P g = -3m$ for some
positive $m$.  Then there is a function $h$ on $E$ that has poles only
at $P$ and at $\infty$ such that $\ord_P (g - h^3 + fh) > -3m$.
Replacing $z$ with $z-h$ and $g$ with $g - h^3 + fh$, we find that
we have reduced the order of the pole of $g$ at $P$.  Repeating this
process, we find that we may assume that for every finite pole $P$ of $g$,
the order of $g$ at $P$ is not a multiple of~$3$.

Suppose $g$ has a pole at $\infty$, and suppose $\ord_\infty g$ is
less than $-3n/2$.  If $\ord_\infty g$ is a multiple of $3$ and is
less than $-3$, then we can find a function $h$, with poles only
at~$\infty$, such that $\ord_\infty (g - h^3 + fh) > \ord_\infty g$.
Again we may replace $z$ with $z-h$ and $g$ with $g - h^3 + fh$ to
reduce the order of the pole of $g$ at $\infty$.  Repeating this
procedure, we may assume that if $\ord_\infty g$ is less than
$-3n/2$ and less than $-3$, then $\ord_\infty g$ is not a multiple of~$3$.

Let us say a pair $(f,g)$ of functions on $E$ is 
\emph{well-conditioned} at a point $P$ of $E$ if one of the following
conditions holds: either
\begin{enumerate}
\item the order $\ord_P g$ of $g$ at $P$ is not a multiple of $3$, or
\item we have $2\ord_P g \ge 3\ord_P f$.
\end{enumerate}
We have shown that every triple cover of $E$
has a model $z^3 - fz = g$ such that $f$ has no poles outside $\infty$
and no multiple zeros anywhere, and such that $(f,g)$ is 
well-conditioned at every finite pole of $g$.  Furthermore,
the model can be made to satisfy the additional requirement that 
$(f,g)$ be well-conditioned at $\infty$, unless $f$ is constant 
and $g$ has a triple pole at $\infty$.

\subsection{Contributions to the different}
\label{SS-different}
Suppose $L/K$ is a field extension of the type considered above,
given in the standard form $z^3 - fz - g = 0$ described in the
preceding section.  
Given a point $P$ on $E$, we would
like to calculate the contribution at $P$ to the different of
the extension $L/K$.
The basic fact we will use is that if $L/K$ is a degree-$3$ Artin-Schreier
extension of local fields given by an equation $z^3 - z = h$, 
where the valuation of $h$ is $n$, then the
degree of the different is zero if $n\ge 0$ and is
$2 - 2n$ if $n<0$ and $n\not\equiv 0 \bmod 3$.
(This follows from~\cite[Prop.~III.7.10]{Stichtenoth}, for example.)
Since the contribution to the different is stable under base
extension, we may assume that the base field $k$ is the algebraic
closure of~$\F_3$.

Suppose we are given $P$ on $E$ with $\ord_P f$ even. 
Note that by the way we normalized $f$ and $g$, either
$2\ord_P g \ge 3\ord_P f$ or $\ord_P g \not\equiv 0 \bmod 3$,
except in the case when $P = \infty$ and $\ord_P g = -3$
and $\ord_P f = 0$.  
Let us suppose we are {\em not\/} in this exceptional case.
Note that if $f$ is nonconstant we will not be in the exceptional case.

In the completion of $K$ at $P$ the function $f$ is a square, say
$f = s^2$ for some $s \in K_P$.  Locally at $P$ the extension 
$L_P / K_P$ is given by the equation $w^3 - w = g/s^3$, and the 
valuation of $g/s^3$ is $\ord_P g - (3/2)\ord_P f$.
Since we are not in the exceptional case, this valuation is
either positive or is not a multiple of~$3$.
Thus the contribution to the different at $P$ is
$$
\begin{cases}
0 & \text{if $3\ord_P f - 2\ord_P g \le 0$;} \\
2 + 3\ord_P f - 2\ord_P g & \text{if $3\ord_P f - 2\ord_P g > 0$.}
\end{cases}
$$
In particular, note that when $\ord_P f$ is even the
contribution at $P$ 
to the different is even, and is at least $4$ if it is nonzero.

Suppose we are given $P$ on $E$ with $\ord_P f$ odd.
Suppose also that $(f,g)$ is well-conditioned at $P$,
so that either
$2\ord_P g \ge 3\ord_P f$ or $\ord_P g \not\equiv 0 \bmod 3$.
Then our completed extension $L_P/K_P$ fits into a diagram
$$\begin{matrix}
L_P & \longrightarrow & L_P'\\
\Buparrow & & \Buparrow \\
K_P & \longrightarrow & K_P'
\end{matrix}
$$
where $K_P'$ and $L_P'$ are obtained from $K_P$ and $L_P$ by
adjoining a square root $s$ of~$f$.
The extension $L_P'/K_P'$ is given by the equation
$w^3 - w = g/s^3$.  Let $P'$ be the prime of $K_P'$.
If the $P'$-adic valuation of $g/s^3$ is nonnegative
(that is, if $2\ord_P g - 3\ord_P f \ge 0$)
then there is no ramification in $L_P'/K_P'$, and so the ramification
in $L_P/K_P$ is tame.  In this case the contribution at $P$
to the different of $L/K$ is~$1.$
On the other hand,
if $2\ord_P g - 3\ord_P f < 0$ then the $P'$-adic valuation of $g/s^3$
is negative and not a multiple of~$3$, so 
the Galois extension $L_P'/K_P'$ is totally ramified.
In particular, $L_P'$ is a field.  Let $\frakp$ be the prime of $L_P'$
lying over $P$.  Then $P' = \frakp^3$, the prime of
$L_P$ is $\frakp^2$, and $P = \frakp^6$.

Let us calculate the different of $L_P' / K_P$ in two different ways.
First of all, we note that the different of $L_P'/K_P'$ is
$$\delta_{L_P'/K_P'} = \frakp^{2 + 2m}$$
where $m = 3\ord_P f - 2\ord_P g$.  Next, we note that the 
extension $K_P'/K_P$ is tamely ramified, so its different
is $\delta_{K_P'/K_P} = P' = \frakp^3.$
Likewise, the different of $L_P'/L_P$ is $\delta_{L_P'/L_P} = \frakp$.
Thus, if the different of $L_P/K_P$ is $\delta_{L_P/K_P} = (\frakp^2)^n$,
we have 
$$\delta_{L_P/K_P} \delta_{L_P'/L_P} = \delta_{K_P'/K_P}\delta_{L_P'/K_P'},$$
so that
$\frakp^{2n + 1} = \frakp^{3 + 2 + 2m}.$
It follows that the contribution $n$ of the different of $L/K$ at $P$
is $n = 2 + 3\ord_P f - 2\ord_P g$.

Thus, when $\ord_P f$ is odd and $(f,g)$ is well-conditioned at $P$,
the contribution to the different at $P$ is
$$
\begin{cases}
1 & \text{if $3\ord_P f - 2\ord_P g \le 0$;} \\
2 + 3\ord_P f - 2\ord_P g & \text{if $3\ord_P f - 2\ord_P g > 0$.}
\end{cases}
$$

Now suppose $P$ is a point for which $\ord_P f$ is odd and $(f,g)$ is
not well-conditioned at $P$.  Since $(f,g)$ is well-conditioned at
infinity unless $f$ is constant, we know that $P$ is a finite point;
furthermore, our normalization of $f$ and $g$ shows that $\ord_P f = 1$.
Also, since $(f,g)$ is well-conditioned at every finite pole of $g$, 
we see that $\ord_P g$ is non-negative.  The only possibility 
is that $\ord_P g = 0$.

Let $c_P\in \kbar$ be the negative of the cube root of the value of $g$
at~$P$.  Then $g+c_P^3 - c_P f$ has a zero at $P$, and 
$(f,g+c_P^3 - c_P f)$ is well-conditioned at $P$.  Note that over
$\kbar$ the triple cover $z^3 - fz = g+c_P^3 - c_P f$ is isomorphic
to the triple cover $z^3 - fz = g$, so we have reduced
the computation of the differential contribution to the case
we considered earlier.

To summarize: If $\ord_P f$ is odd, let $c_P$ be an element of 
$\kbar$ such that either $2\ord_P (g+c_P^3 - c_P f) \ge 3\ord_P f$ or 
$\ord_P (g+c_P^3 - c_P f) \not\equiv 0\bmod 3$.  
Let $g_P = g+c_P^3 - c_P f$.  Then the contribution to the 
different at $P$ is
$$
\begin{cases}
1 & \text{if $3\ord_P f - 2\ord_P g_P \le 0$;} \\
2 + 3\ord_P f - 2\ord_P g_P & \text{if $3\ord_P f - 2\ord_P g_P > 0$.}
\end{cases}
$$
In particular, when $\ord_P f$ is odd the contribution at $P$ to the
different is odd.

\section{An argument on Hermitian forms}
\label{S-Hermitian}

In this section we prove a theorem of Savitt~\cite{Savitt}:

\begin{thm}
\label{T-savitt}
There is no genus-$4$ curve over $\F_8$ with exactly $27$ rational points. 
\end{thm}

\begin{proof}
Suppose such a curve $C$ existed.  
It has defect $2$ and so 
we know that $C$ must be of type $[m,m,m,m-2]$ or of type 
$$[m + (-1+\sqrt{5})/2, m + (-1-\sqrt{5})/2,
   m + (-1+\sqrt{5})/2, m + (-1-\sqrt{5})/2],$$
where $m = 5$.
Corollary~\ref{C-type2} eliminates the first possibility,
so $C$ must be of the latter type.  It follows that the 
Weil polynomial 
of $C$ must be~$f^2$, where $f = x^4-9x^3+35x^2-72x+64$. 
Our proof of Savitt's theorem is completed by the following 
proposition, which shows that every principal polarization 
of an abelian variety with Weil polynomial $f^2$ 
is decomposable. 
\end{proof}

\begin{prop}
\label{P-decomposable}
There is exactly one abelian variety $A$ over $\F_8$ with Weil
polynomial $f$. Up to isomorphism, the variety $A$ has exactly 
one principal polarization $\lambda$. Furthermore, up to 
isomorphism there is exactly one principally polarized abelian 
variety over $\F_8$ with Weil polynomial~$f^2$, and it is 
isomorphic to $(A\times A, \lambda\times\lambda)$. 
\end{prop}

Let $K$ be the quartic number field defined by the polynomial $f$ 
and let $\calO_K$ denote the ring of integers of $K$.  Our proof 
of Proposition~\ref{P-decomposable} will depend on a result about 
Hermitian forms over $\calO_K$.  We will state this result now and 
use it in the proof of Proposition~\ref{P-decomposable}, but 
we will postpone its proof until later in this section.

We will see that $K$ is the totally imaginary biquadratic
extension $\Q(\sqrt{-3},\sqrt{5})$ of the totally real field 
$\Kp = \Q(\sqrt{5})$;  we refer to the nontrivial 
automorphism of $K$ over $\Kp$ as {\em complex conjugation\/}, 
and we denote the complex conjugate of $x\in K$ by $\overline{x}$. 
Let $M_2(\calO_K)$ denote the ring of $2$-by-$2$ matrices over $\calO_K$. 
If $C$ is an element of $M_2(\calO_K)$ we let $C^*$ denote its 
conjugate-transpose.

\begin{prop}
\label{P-reduction}
Suppose $A$ is an invertible Hermitian matrix in $M_2(\calO_K)$ that 
is totally positive {\rm(}meaning that all of the roots of its minimal 
polynomial are totally positive algebraic numbers\/{\rm)}.  Then there 
is an invertible $C\in M_2(\calO_K)$ such that $A = C^*C$. 
\end{prop}

Let us assume this result for the time being, and proceed 
with the proof of Proposition~\ref{P-decomposable}.

\begin{proof}[Proof of Proposition~{\rm\ref{P-decomposable}}]
We begin by setting some notation related to the number field~$K$.

Let $\pi$ be a root of $f$ in~$K$, let $\pibar = 8/\pi$, 
and let $R = \Z[\pi,\pibar]$. Let $\varphi = \pi + \pibar - 4$ 
and let $\zeta = 17 - 6\pi + \pi^2 - 3\pibar$. It is easy to 
check that then $\varphi^2 - \varphi - 1 = 0 $ and 
$\zeta^2 + \zeta + 1 = 0$, and from these relations 
we see that $K$ is isomorphic to $\Q(\sqrt{-3},\sqrt{5})$
and that $R$ is the full ring of integers of $K$.

It is not hard to show that the Dedekind domain $R$ is a PID~\cite{Lakein}; 
in fact, Lemma~\ref{L-EuclideanAlg} below shows that $R$ is norm-Euclidean.

Note that the middle coefficient of $f$ is coprime to $8$, so 
there is an isogeny class of ordinary abelian varieties over $\F_8$ 
with Weil polynomial~$f$.  In fact, according to a result
of Deligne~\cite{Deligne}, the abelian varieties in this isogeny 
class correspond to the isomorphism classes of $R$-modules that 
can be embedded in $K$ as lattices. Since $R$ is the full ring 
of integers of $K$ and since $R$ has class number~$1$, there is 
exactly one such isomorphism class of $R$-modules, and therefore 
there is exactly one abelian variety $A$ over $\F_8$ with 
Weil polynomial~$f$. This proves the first statement
of the proposition.

Theorem~1.3 of~\cite{Howe} shows that the abelian variety $A$ 
has a principal polarization~$\lambda$.  Now suppose $\mu$ 
is another principal polarization of $A$.  Then we know 
from~\cite[Application~III, pp.~208--210]{Mumford} 
(see especially the final paragraph) that there is a totally 
positive unit $u$ of the maximal real subfield $K^+$ of $K$ 
such that $\mu = \lambda u$.  But every totally positive unit 
of $K^+ = \Q(\sqrt{5})$ is an even power of the
fundamental unit $\varphi$,  so there is a unit $v$ 
of $K^+$ with $u = v^2 = v\overline{v}$. Then the 
automorphism $v$ of $A$ gives an isomorphism of the 
polarized varieties $(A,\lambda)$ and $(A,\mu)$.  
This proves the second statement of the proposition.

Applying Deligne's theorem again, we find that the abelian 
varieties over $\F_8$ with Weil polynomial $f^2$ 
correspond to the isomorphism classes of $R$-modules that 
can be embedded as lattices in the $K$-vector space~$K\times K$. 
Since $R$ is a Dedekind domain, such modules are determined 
up to isomorphism by their Steinitz classes in the class 
group of~$R$.  But the class group of $R$ is trivial, so there 
is only one such $R$-module.  Thus, the only abelian variety 
with Weil polynomial $f^2$ is $A\times A$.  

Now suppose that $\mu$ is a principal polarization on $A\times A$. 
Let $\alpha$ be the automorphism 
$\mu^{-1}\circ(\lambda\times\lambda)$ of~$A\times A$. 
Using the results of~\cite[Application~III, pp.~208--210]{Mumford} 
again, we see that $\alpha$ is fixed by the Rosati involution 
associated to $\lambda\times\lambda$ and that $\alpha$ is totally 
positive, meaning that all of the roots (in the algebraic closure 
of $K$) of the minimal polynomial of $\alpha$ are totally positive
algebraic numbers.  If we identify $\End(A\times A)$ with the 
ring $M_2(R)$ of $2$-by-$2$ matrices over~$R$ in the obvious way, 
then the Rosati involution is the conjugate-transpose involution, 
so we see that $\alpha$ is identified with a totally positive 
Hermitian matrix $A$ of determinant~$1$.   Thus, to show that 
$\mu$ is isomorphic to $\lambda\times\lambda$, we must show 
that every such Hermitian matrix can be written $C^*C$, 
where $C\in M_2(R)$ is nonsingular and where $C^*$ is the 
conjugate transpose of~$C$. But this is exactly the statement 
of Proposition~\ref{P-reduction}. 
\end{proof}

Before we prove Proposition~\ref{P-reduction} we must set some notation 
and give a Euclidean algorithm for the ring~$R$.
 
Let $L$ be the subfield $\Q(\zeta)$ of $K$ and let $\calO_L$ be 
the ring of integers $\Z[\zeta]$ of $L$. Let $\phi$ be the real 
number $(1 + \sqrt{5})/2$ and let $\psi_1$ and $\psi_2$ be two 
distinct embeddings of $K$ into $\C$ that are not complex conjugates
of one another.  If $z$ is a complex number, we let $|z|$  be 
its magnitude and we let $||z||$ be its norm, so that  
$||z|| = |z|^2 = z \bar{z}$.

\begin{lem}
\label{L-Euclidean}
For every $x$ in $K$ there is a $y$ in $R$ such that 
$$\Norm_{K/\Q} (x - y) \le 5/9$$ 
and such that 
$$||\psi_i(x - y)|| \le \phi^4/3 \quad\text{for  $i = 1,2.$}$$ 
\end{lem}

\begin{proof}
Let $D$ be the set of elements of $L$ whose norm to $\Q$ is at most~$1/3$. 
Then for every $x$ in $L$ there is a $y$ in $\calO_L$ such that $x - y$ 
lies in $D$. (It is easiest to see this by embedding $L$ in the complex 
numbers, so that $D$ becomes the intersection of $L$ with the disk at the 
origin of radius $1/\sqrt{3}$.  The latter disk clearly contains a 
fundamental region for the lattice $\calO_L$.)

Write $x = x_1 + x_2 \varphi$ for $x_1, x_2\in L$.  Choose $y_1$ 
and $y_2$ in $O_L$ such that $x_1-y_1$ and $x_2-y_2$ lie in $D$. 
Let $y = y_1 + y_2 \varphi$ and let $z_1 = x_1-y_1$ and $z_2 = x_2-y_2$. 
Then 
\begin{align*} 
\Norm_{K/\Q}(x - y)  & = \Norm_{L/\Q}( \Norm_{K/L}( z_1 + z_2 \varphi)) \\
                     & = \Norm_{L/\Q}(z_1^2 + z_1z_2 - z_2^2) \\
                     & = || \psi_1(z_1^2 + z_1z_2 - z_2^2) ||  \\ 
                     & = || d_1^2 + d_1d_2 - d_2^2 || 
\end{align*} 
where $d_1 = \psi_1(z_1)$ and $d_2 = \psi_1(z_2)$ are complex numbers 
that lie in the disk around the origin of radius $1/\sqrt{3}$. 
But an easy maximization argument shows that the maximum value 
of $|| d_1^2 + d_1 d_2 - d_2^2 ||$ for $d_1$, $d_2$ in this disk is $5/9$.

Also, $\psi_i(x - y)$  is equal to either $d_1 + d_2 (1 + \sqrt{5})/2$ 
or $d_1 + d_2 (1 - \sqrt{5})/2$, depending on the image of $\varphi$ 
under~$\psi_i$.  Since $|d_1|$ and $|d_2|$ are both at most $1/\sqrt{3}$, 
we see that 
$$  | \psi_i(x - y) |  \le \frac{1}{\sqrt{3}} \frac{3 + \sqrt{5}}{2} 
                         = \frac{\phi^2}{\sqrt{3}} $$
so that  $|| \psi_i(x - y) || \le \phi^4/3$.
\end{proof}

We note that the $5/9$ in Lemma~\ref{L-Euclidean} could be reduced 
to $4/9$ if we used a hexagonal fundamental domain for the lattice 
$\calO_L$ in place of the disk $D$, but doing so takes some effort 
and does not help much in the end.

\begin{lem}
\label{L-EuclideanAlg}
Suppose $n$ and $d$ are elements of $\calO_K$, with $d$ nonzero.  
Then there are elements $q$ and $r$ of $O_K$ such that $n = qd + r$ 
and such that 
$$\Norm_{K/Q} (r) \le (5/9)\Norm_{K/Q} (d)$$ 
and 
$$||\psi_i(r)|| \le  \frac{\phi^4}{3} ||\psi_i(d)||$$ for $i = 1,2$. 
\end{lem}

\begin{proof}
Apply Lemma~\ref{L-Euclidean} to $x=n/d$, and let $q$ be the 
resulting $y$. Then let $r = n - qd$.  The lemma follows from 
the inequalities of Lemma~\ref{L-Euclidean}. 
\end{proof}

We are now ready to prove Proposition~{\rm\ref{P-reduction}}.

\begin{proof}[Proof of Proposition~{\rm\ref{P-reduction}}]
Write 
$$A= \left[\ \begin{matrix}
\alpha & \betabar\\
\beta  & \gamma
\end{matrix}\ \right]$$
for $\alpha$, $\gamma$ in the ring of integers of the maximal real 
subfield $\Kp=\Q(\varphi)$ of K and for $\beta$ in $R$. Our strategy
will be to modify $A$ by invertible matrices $C$ (that is, to 
replace $A$ with $C^*AC$) in order to make the norm of the upper 
left hand element of $A$ as small as possible.

The determinant of $A$ is a totally positive unit in $\Kp$, and so is 
an even power of the fundamental unit $\varphi$. By modifying $A$ 
by a matrix $C$ of the form 
$$\left[\ \begin{matrix}
\varphi^i & 0 \\
0 & 1
\end{matrix}\ \right]$$
we may assume that $A$ has determinant~$1$. Then by modifying 
$A$ by a power of the matrix 
$$\left[\ \begin{matrix}
\varphi & 0\\
0 & \varphi^{-1}
\end{matrix}\ \right]$$
we can ensure that the element $\alpha$ of $\calO_{\Kp}$ has the 
property that 
$$\frac{1}{\phi^2} \le \frac{\psi_1(\alpha) }{ \psi_2(\alpha)}  
                   \le  \phi^2 .$$
Another way of expressing this is to say that 
\begin{equation} 
\label{EQ-alphabound} 
\frac{1}{\phi^2} \le \frac{\psi_i(\alpha)^2 }{\Norm_{\Kp/\Q}(\alpha)}  
                 \le \phi^2 \quad\text{for $i = 1, 2$.} 
\end{equation}

Apply Lemma~\ref{L-EuclideanAlg} with $n = \beta$ and $d = \alpha$ to get 
a $q$ and an~$r$ with $||\psi_i(r)|| \le (\phi^4/3)||\psi_i(\alpha)||$
and with $\Norm_{K/\Q}(r) \le (5/9) \Norm_{K/\Q}(\alpha)$.
If we set 
$$C= 
\left[\ \begin{matrix}
1 &  -\qbar\\
0 &  1
\end{matrix}\ \right]$$
then 
$$C^* A C  = 
\left[\ \begin{matrix}
\alpha & \rbar\\
r      & \gamma'
\end{matrix}\ \right]$$
for some $\gamma'$ in $\calO_\Kp$.   Replace $\beta$ with $r$ and 
$\gamma$ with~$\gamma'$, so that now we have
\begin{equation}
\label{EQ-betabound}
||\psi_i(\beta)|| \le \frac{\phi^4}{3} ||\psi_i(\alpha)||
                      \quad\text{for $i=1,2$}
\end{equation}
and 
\begin{equation}
\label{EQ-betabound2}
\Norm_{K/\Q}(\beta) \le (5/9) \Norm_{K/\Q}(\alpha).
\end{equation}

Let $B = \beta \betabar$, so that $B$ is an element of $\calO_\Kp$. 
Note that we have $\alpha\gamma - B = 1$, so
$$  \psi_i(\alpha) \psi_i(\gamma) = 1 + \psi_i(B) \quad \text{for $i = 1,2 $}$$ 
and 
\begin{equation} 
\label{EQ-gammabound}
\psi_i(\gamma) / \psi_i(\alpha) = 
   1/\psi_i(\alpha)^2 + \psi_i(B)/\psi_i(\alpha)^2 \quad \text{for $i = 1,2$.} 
\end{equation}

Now  let 
\begin{align*}
b_1 &= \psi_1(B) / \psi_1(\alpha^2) \\
b_2 &= \psi_2(B) / \psi_2(\alpha^2) \\
c_1 &= 1/\psi_1(\alpha^2)\\
c_2 &= 1/\psi_2(\alpha^2)
\end{align*}
so that equation~(\ref{EQ-gammabound}) becomes
$$\psi_1(\gamma) / \psi_1(\alpha)  = b_1 + c_1 \qquad\text{and}\qquad
  \psi_2(\gamma) / \psi_2(\alpha)  = b_2 + c_2.$$
Multiplying these last two equalities gives
\begin{equation}
\label{EQ-gammaalphabound}
\Norm_{\Kp/\Q}(\gamma/\alpha) = b_1 b_2 + b_1 c_2 + b_2 c_1 + c_1 c_2. 
\end{equation} 
Note that 
\begin{equation} 
\label{EQ-bbbound}
  b_1 b_2 = \frac{\Norm_{\Kp/\Q}(B)}  {\Norm_{\Kp/\Q}(\alpha)^2}
          = \frac{\Norm_{K/\Q}(\beta)}{\Norm_{K/\Q}(\alpha)} 
        \le \frac{5}{9}
\end{equation}
(where the final inequality comes from~(\ref{EQ-betabound2}))
and 
\begin{equation}
\label{EQ-ccbound}
  c_1 c_2 = \frac{1}{\Norm_{\Kp/\Q}(\alpha^2)}.
\end{equation}
Furthermore, from inequality~(\ref{EQ-alphabound}) we see that 
\begin{equation} 
\label{EQ-cbound} 
c_1  \le \frac{\phi^2}{\Norm_{\Kp/\Q}(\alpha)} \qquad\text{and}\qquad 
c_2  \le \frac{\phi^2}{\Norm_{\Kp/\Q}(\alpha)},
\end{equation}
and from inequality~(\ref{EQ-betabound}) we see that 
\begin{equation} 
\label{EQ-bbound} 
b_1 = \frac{||\psi_1(\beta)||}{||\psi_1(\alpha)||} \le \frac{\phi^4}{3} 
\qquad\text{and}\qquad
b_2 = \frac{||\psi_2(\beta)||}{||\psi_2(\alpha)||} \le \frac{\phi^4}{3}. 
\end{equation}

If we view $b_1$, $b_2$, $c_1$, and $c_2$ as non-negative real 
variables subject only to the conditions expressed in 
equations~(\ref{EQ-bbbound}), (\ref{EQ-ccbound}), (\ref{EQ-cbound}), 
and~(\ref{EQ-bbound}), and if we maximize $b_1c_1 + b_2c_2$ subject 
to these conditions, we find that the maximum value occurs when 
$b_1 = \phi^4/3$ and $c_1 = \phi^2/\Norm_{\Kp/\Q}(\alpha)$. 
Thus we have 
\begin{equation} 
\label{EQ-crosstermbound}
b_1c_1 + b_2c_2 \le \frac{\phi^4}{3} \frac{\phi^2}{\Norm_{\Kp/\Q}(\alpha)}
       + \frac{(5/9)}{(\phi^4/3)} \frac{(1/\phi^2)}{\Norm_{\Kp/\Q}(\alpha)}
         \le \frac{6.08}{\Norm_{\Kp/\Q}(\alpha)}. 
\end{equation} 
Let $\epsilon = 1/\Norm_{\Kp/\Q}(\alpha)$.  Then by combining 
the relations~(\ref{EQ-gammaalphabound}), (\ref{EQ-bbbound}),
(\ref{EQ-ccbound}) and~(\ref{EQ-crosstermbound}) we find that 
$$\Norm_{\Kp/\Q}(\gamma/\alpha) \le \epsilon^2  + 6.08 \epsilon + 5/9.$$

If $\Norm_{\Kp/\Q}(\alpha) \ge 15$ then $\epsilon < 0.07$ and 
$\Norm_{\Kp/\Q}(\gamma/\alpha) < 1$.  Then we can modify $A$ by 
$$\left[\ \begin{matrix}
0 & 1 \\
1 & 0 \end{matrix}\ \right]$$
to exchange $\alpha$ and $\gamma$, and this decreases the norm of the 
upper left hand element of $A$.

We repeat this procedure until we reach the point where
$\Norm_{\Kp/\Q}(\alpha) \le 14$.  Up to Galois conjugacy
there are $5$ possible values 
for $\alpha$: namely, $1$, $2$, $2+\varphi$, $3$, and $3+\varphi$. 
Inserting the appropriate values of $c_1$ and $c_2$ into 
equation~(\ref{EQ-gammaalphabound}) and maximizing over $b_1$ and $b_2$, 
we find that we once again get $\Norm_{\Kp/\Q}(\gamma/\alpha) < 1 $ 
except when $\alpha = 1$ or $\alpha = 2$ or $\alpha = 2 + \varphi$.

Note that $1 + \beta \betabar$ is divisible by~$\alpha$.  If 
$\alpha = 2 + \phi$ then there are $6$ possible residue classes modulo 
$\alpha$ that $\beta$ could lie in.  By modifying $A$ by a power of 
$$\left[\ \begin{matrix}
1 & 0 \\
0 & \zeta\end{matrix}\ \right]$$
we can arrange that  $\beta = 2$.  Then $A$ must be the matrix 
$$\left[\ \begin{matrix} 
2 + \varphi & 2 \\
2           & 2 +\overline{\varphi}
\end{matrix}\ \right].$$
Modifying $A$ by 
$$\left[\ \begin{matrix}
1-\varphi & 1 \\
\varphi   & -1
\end{matrix}\ \right]$$ 
gives us the identity.

For $\alpha = 2$ there are three possible residue classes for $\beta$. 
By modifying $A$ by a power of 
$$\left[\ \begin{matrix}
1 & 0 \\ 
0 & \zeta
\end{matrix}\ \right]$$
we can arrange for  $\beta$ to be $1$.  Then $\gamma = 1$, 
so again we can reduce the norm of the upper left hand corner 
of $A$ by interchanging $\alpha$ and $\gamma$.

We finally get to the case $\alpha = 1$.  If $\alpha = 1$ then 
we can reduce $\beta$ to be $0$.  Then we find $\gamma = 1$,
so that we have reduced $A$ to the identity matrix.
\end{proof}

\section*{Appendix I: Sample Magma output}
In this appendix we reproduce the output produced by our Magma program
for the case $q=4$, $g=5$, $N=18$.

\begin{verbatim}
Magma V2.9-10     Fri Sep 27 2002 17:22:30    [Seed = 729997397]
Type ? for help.  Type <Ctrl>-D to quit.
> load "CheckQGN.magma";
Loading "CheckQGN.magma"
Loading "DeficientPolynomialList.magma"
> CheckQGN(4,5,18);
[ 18, 0, 4, 81, 164 ]
[
    <x + 3, 3>,
    <x^2 + 4*x + 2, 1>
]
ELIMINATED: Not Weil polynomial.

[ 18, 0, 5, 74, 187 ]
[
    <x + 3, 2>,
    <x^3 + 7*x^2 + 14*x + 7, 1>
]
ELIMINATED: resultant=1 method.
Splitting = [ 1 ]

[ 18, 0, 6, 67, 210 ]
[
    <x + 1, 1>,
    <x + 2, 1>,
    <x + 3, 2>,
    <x + 4, 1>
]
ELIMINATED: resultant=2 method.
Splitting = [ 2 ]
Reasons: point counts, Riemann-Hurwitz

[ 18, 0, 6, 68, 200 ]
[
    <x + 3, 1>,
    <x^2 + 5*x + 5, 2>
]
ELIMINATED: resultant=1 method.
Splitting = [ 1 ]

[ 18, 0, 7, 60, 232 ]
[
    <x + 2, 2>,
    <x + 4, 1>,
    <x^2 + 5*x + 5, 1>
]
ELIMINATED: resultant=1 method.
Splitting = [ 1, 2 ]

[ 18, 1, 0, 86, 168 ]
[
    <x + 1, 1>,
    <x + 3, 4>
]
ELIMINATED: resultant=2 method.
Splitting = [ 1 ]
Reasons: point counts, Riemann-Hurwitz

[ 18, 1, 1, 79, 190 ]
[
    <x + 2, 1>,
    <x + 3, 2>,
    <x^2 + 5*x + 5, 1>
]
ELIMINATED: resultant=1 method.
Splitting = [ 1 ]

[ 18, 1, 2, 71, 222 ]
[
    <x + 2, 3>,
    <x + 3, 1>,
    <x + 4, 1>
]
ELIMINATED: resultant=1 method.
Splitting = [ 2 ]

> quit;

Total time: 4.940 seconds, Total memory usage: 9.22MB
\end{verbatim}

\section*{Appendix II: New arguments for two cases}

\noindent
{\bf A new argument for the case $q=3$, $g=6$, $N=15$.}
Suppose that $C$ is a genus-$6$ curve over $\F_3$ with exactly $15$
rational points. In Section~\ref{SS-3} we showed that the real 
Weil polynomial of the Jacobian $J$ of $C$  must be 
$$h = (x+2)^2 (x+3) (x^3 + 4x^2+x-3).$$
Let $F$ be the unique elliptic curve over $\F_3$ with real Weil 
polynomial $x+2$, and let $B$ be the abelian surface $F\times F$.  
Note that $B$ is the only abelian surface over $\F_3$ with real
Weil polynomial $(x+2)^2$.  We can show that there is an injection 
$B\hookrightarrow J$ such that the canonical polarization on $J$ pulls
back to a polarization $\mu$ on $B$ whose degree is $9$.  By looking at
the degree-$9$ polarizations of $B$, we see that there will be an 
injection $F\hookrightarrow B$ such that the pullback of $\mu$ to $F$
is a polarization $\lambda$ of degree $1$ or $4$.  Now consider the 
composition 
$$F\hookrightarrow B \hookrightarrow J.$$
The canonical polarization on $J$ pulls back via this composition to
the polarization $\lambda$ of $F$, and it follows that there is a map 
from $C$ to $F$ of degree $1$ or $2$.  Certainly there is no such map 
of degree $1$.  But there are no such maps of degree $2$ either, 
because $F$ has $6$ rational points and $C$ is supposed to have $15$.
Therefore there is no genus $6$ curve over $\F_3$ with real Weil 
polynomial equal to $h$.

\medskip
\noindent
{\bf A new argument for the case $q=27$, $g=4$, $N=65$.}
The appendix to~\cite{Savitt} shows that a genus-$4$ curve over 
$\F_{27}$ with $65$ rational points must have real Weil polynomial
$(x+7)(x+10)^3$.  Suppose $C$ is such a curve, and let $J$ be its 
Jacobian.  Let $F$ be the unique elliptic curve over $\F_{27}$ with 
real Weil polynomial $x+10$.  Note that $F^3$ is the unique abelian 
threefold over $\F_{27}$ with real Weil polynomial $(x+10)^3$.  We see
from Lemma~\ref{L-resultant} that there is a degree-$9$ isogeny 
$E \times F^3 \to J$, and that the pullback to $F^3$ of the canonical
principal polarization on $J$ is a degree-$9$ polarization on $F^3$. 
Using the knowledge of the isomorphism classes of principal 
polarizations of $F^3$ that we obtain from Hoffmann's classification
\cite{Hoffmann} of the rank-$3$ unimodular lattices over 
$\Z[\sqrt{-2}]$, together with an easy argument that shows that every
degree-$9$ polarization of $F^3$ is the pullback of a principal 
polarization on $F^3$ via a $3$-isogeny, we can write down 
representatives for all of the isomorphism classes of degree-$9$
polarizations on $F^3$.  For each representative $\mu$, we check that 
there is an embedding $F\hookrightarrow F^3$ such that the pullback of
$\mu$ to $F$ is a polarization of degree $1$ or $4$.  It follows that 
there must be a map from $C$ to $F$ of degree $1$ or $2$.  A degree-$1$
map would be impossible, so there must be a degree-$2$ map.  But it is
not hard to adapt the method explained in Section~\ref{SS-27.4.66} to
enumerate the genus-$4$ double covers of $F$, and to verify that none
of them has $65$ points.  Therefore there is no genus-$4$ curve over
$\F_{27}$ having $65$ points.


\end{document}